\documentclass[10pt,reqno]{amsart}
\usepackage[utf8]{inputenc}
\usepackage[T1]{fontenc}
\usepackage[english]{babel}
\usepackage[margin=2.5cm]{geometry}
\usepackage{amsmath, amsfonts, amsthm, amssymb, graphicx}
\usepackage[usenames,dvipsnames,svgnames,table]{xcolor}
\usepackage[colorlinks=true, hyperindex, pagebackref=true, citecolor=magenta,linkcolor=Blue]{hyperref}
\usepackage[all]{xy}
\usepackage{enumitem}
\setlist[enumerate]{label=\rm(\arabic*)}
\newtheorem{lemma}[subsection]{Lemma}
\newtheorem{proposition}[subsection]{Proposition}
\newtheorem{theorem}[subsection]{Theorem}
\newtheorem{remark}[subsection]{Remark}
\newtheorem{definition}[subsection]{Definition}

\newtheorem{corollary}[subsection]{Corollary}

\newcommand{\bbR}{\mathbb{R}}
\newcommand{\bbZ}{\mathbb{Z}}

\newcommand{\bbN}{\mathbb{N}}
\newcommand{\bfF}{\mathbf{F}}

\newcommand{\Sid}{\mathcal{S}_{\id}(f)}

\newcommand{\F}{\mathcal{F}(M,P)}

\newcommand{\id}{\mathrm{id}}
\newcommand{\Crit}{\Sigma_f}
\newcommand{\CritP}{\Crit^P}
\newcommand{\CritC}{\Crit^C}

\renewcommand{\ker}{\mathrm{Ker}}

\newcommand{\textr}[1]{\textcolor{black}{#1}}

\title[Homotopy type of stabilizers of smooth functions]{Homotopy type of stabilizers of smooth functions with non-isolated singularities on surfaces}
\author{Bohdan Feshchenko}
\address{Department of algebra and topology, Institute of Mathematics of National Academy of Science of Ukraine,
	Tereshchenkivska, 3, Kyiv, 01601, Ukraine}
\email{fb@imath.kiev.ua}
\date{\today}
\begin{document}
\maketitle
\begin{abstract}
The paper is devoted to the study of homotopy properties of stabilizers of smooth functions on orientable surfaces, i.e., groups of diffeomorphisms of surfaces preserving a given function. For some class of smooth functions which is a generalization of the class of Morse-Bott functions on orientable surfaces, the homotopy type of the path component of the identity map of the stabilizer is completely described.
\end{abstract}
\section{Introduction}
The stabilizers and orbits of smooth functions on compact surfaces are special spaces of smooth maps which arise from the action of the group of diffeomorphisms of surfaces on the space of smooth functions. These spaces naturally appear  in problems of smooth equivalences and  deformations for smooth functions on surfaces. A systematic study of their homotopy properties  has been started  by S.~Maksymenko \cite{Maksymenko:AGAG:2006}.
 We will give precise definitions of stabilizers and orbits below.

Let $M$ be a smooth, compact surface, and $P$ be either a real line $\bbR$ or a circle $S^1$. The group  of diffeomorphisms $\mathcal{D}(M)$ of $M$ acts on the space of smooth $P$-valued functions $C^{\infty}(M,P)$ by the following rule:
$$
\xi: C^{\infty}(M,P)\times \mathcal{D}(M)\to C^{\infty}(M,P),\qquad \xi(f, h) = f\circ h.
$$
For a smooth function $f\in C^{\infty}(M,P)$,  we denote by 
$$
\mathcal{S}(f) = \{h\in \mathcal{D}(M)\,| f\circ h = f \},\quad \mathcal{O}(f) = \{ f\circ h\,|\, h\in \mathcal{D}(M) \}
$$
the {\it stabilizer} and the {\it orbit} of $f$ with respect to the action $\xi$. Endow strong Whitney topologies on $\mathcal{D}(M)$ and $C^{\infty}(M,P)$; these topologies induce some topologies on $\mathcal{S}(f)$ and $\mathcal{O}(f)$. We also denote by $\mathcal{D}_{\id}(M)$ and $\mathcal{S}_{\id}(f)$ path components  of $\mathcal{D}(M)$ and $\mathcal{S}(f)$ containing $\id_M$, and by $\mathcal{O}_f(f)$ the path component of $\mathcal{O}(f)$ containing $f$.

For   $P$-valued Morse functions on compact surfaces the homotopy types of $\Sid$ and $\mathcal{O}_f(f)$ were described by S.~Maksymenko \cite{Maksymenko:AGAG:2006} and E.~Kudryavtseva \cite{Kudryavtseva:MathNotes:2012, Kudryavtseva:MatSb:2013}. 
In particular,
 S.~Maksymenko showed that $\Sid$ is contractible if  $f$ has either  at least one saddle point or $M$ is non-orientable; otherwise, $\Sid$ is homotopy equivalent to $S^1$.
It was proved that the orbit $\mathcal{O}_f(f)$ is homotopy equivalent to $m$-torus $T^m$ if $M$ is aspherical, to $S^2$ if $M = S^2$ and $f$ has exactly $2$ critical points, and to $\mathrm{SO}(3)\times T^m$ otherwise, for some $m\geq 0$  depending on $f$. 
Algebraic structures of homotopy groups of orbits and other groups which ``partially'' controls the homotopy type of orbits are well understood for such smooth functions on all compact orientable surfaces except $S^2$ and  remains less understood for functions on certain non-orientable surfaces, such as the Klein bottle and the projective plane.
These results were also generalized to a larger class of smooth functions with isolated singularities;
 more details  can be found in  \cite{Maksymenko:2021:review}.
We  also note that recently homotopy properties of orbits were applied to some questions on persistent homology of Morse functions  by J.~Leygonie and D.~Beers \cite{LeygonieBeers:JACT:2023}.

Our  goal is to generalize results on the homotopy type of stabilizers to more general class
of circle-valued functions on surfaces whose critical points can be non-isolated. The natural class of such functions to consider are Morse-Bott functions.

Morse-Bott functions are the generalization of Morse functions by significantly weakening the conditions on the set of critical points -- a critical set of a Morse-Bott function is a union of submanifolds and each connected component of such union is  ``non-degenerate in the normal direction''. General information about them in the context of Morse theory can be found \cite{Nicolaescu:MorseTh:2011}. 
Such functions are more flexible than Morse functions and can ``capture'' additional information about the symmetries that a manifold may possess, which is especially important for higher-dimensional manifolds.
Morse–Bott functions are  very popular objects of study and have been investigated by many specialists in various problems of topology and its applications.
In particular,
topological classification of Morse-Bott functions on orientable surfaces were studied by  E.~B.~Batista, J.~C.~F.~Costa and I.~S.~Meza-Sarmiento \cite{BatistaCostaMeza:2022} and by J.~Mart\'{\i}nez-Alfaro, I.~S.~Meza-Sarmiento  and R.~Oliveira  \cite{MartinezMezaOliveria:2016}. I.~Gelbukh classified Morse-Bott functions on manifolds with the only $2$  critical values
\cite{Gelbukh:Czechoslovak:2021}. Topological properties foliations with Morse-Bott singularities of codimension-$1$ were investigated  by B.~Sc\'{a}rdua and J.~Seade 
\cite{ScarduaSeade:DGeom:2009,ScarduaSeade:Topol:2011}, and the homotopy properties of diffeomorphisms preserving Morse-Bott foliations on lens spaces were studied by S.~Maksymenko \cite{Maksymenko:2023:MBFol2:2023,MaksymenkoKhokhliuk:2022:MBFol:2022}.

In this paper we consider the class $\F$ of $P$-valued functions\footnote[2]{In the following, we will omit the term ``$P$-valued'' and indicate it only when necessary. Thus, maps from $C^{\infty}(M,P)$ will simply be called functions.} on surfaces with ``generalized Morse-Bott'' singularities  
and describe the homotopy type of $\Sid$ for functions on orientable surfaces from $\F$ (see Theorem \ref{thm:main-theorem} below).  Denote by $C_{\partial}^{\infty}(M,P)$ a subclass of $C^{\infty}(M,P)$ of smooth functions which are locally constant on the boundary $\partial M.$
\begin{definition}\label{def:class-F} 
	 \normalfont A smooth function $f\in C_{\partial}^{\infty}(M,P)$ on $M$ belongs to the class $\mathcal{F}(M, P)$ if it satisfies the following conditions:
	\begin{enumerate}
		\item the set of critical points $\Sigma_f$ of $f$ is a disjoint union of smooth submanifolds of $M$ and $\Sigma_f\subset \mathrm{Int}(M),$
		\item  for each connected component $C$ of $\Sigma_f$ and any critical point $p\in C$,
		the germ $(f,p)$ of $f$ at $p$ is smoothly equivalent 
		\begin{enumerate}
			\item[(a)] to either \textr{the germ at $0\in \bbR^2$ of a homogeneous polynomial} $f_p:\bbR^2\to \bbR$ without multiple factors  with $\deg f_p\geq 2$,
			\item[(b)] or to \textr{the germ at $0\in \bbR^2$ } of  $f_C(x,y) = \pm y^{n_C}$  for some $n_C\in \bbN_{\geq 2}$ depending of $C.$  
		\end{enumerate}
	\end{enumerate}
\end{definition}
For a function $f\in \F$, a connected component $C$ of $\Crit$ is either an isolated critical point \textr{due to (2.a)} or a critical circle \textr{due to (2.b)}. 
{An isolated critical point  $p$ will be called a {\it saddle}, if a polynomial $f_p$ from (2.a) has at least one linear factor. A critical circle $C$ of $f$ can be  extremal ($n_C$ is even) or non-extremal ($n_C$ is odd).}
A local structure of level-sets  of functions from the class $\F$ near their singularities will be discussed  in \textsection \ref{ssec:class-F}. 
Clearly, $\F$ contains the class of Morse-Bott functions, as well as  Morse functions locally constant on $\partial M$.

The following theorem is our main result.
\begin{theorem}\label{thm:main-theorem}
	Let $M$ be a smooth, compact, connected, and \textr{orientable} surface, and let $f$ be a function from $\F$. 
Then $\Sid$ is either contractible or homotopy equivalent to a circle $S^1$. To be more precise, 
$\Sid$ is contractible if either $f$ has at least one  saddle point or $f$ has a \textr{degenerate} isolated extremum; otherwise $\Sid$ is homotopy equivalent to $S^1.$
\end{theorem}
 Thus, for ``almost all'' function from $\F$, the stabilizer $\Sid$ is contractible.
 The combinatorial properties of  functions whose stabilizer $\Sid$ is homotopy equivalent to  $S^1$ are described in the following statement.
\begin{proposition}\label{prop:main}
	Let $M$ be a smooth, compact, connected, and \textr{orientable} surface, and let $f\in \F$ be a function such that $\Sid$ is homotopy equivalent to $S^1$. Then the following  hold:
	\begin{enumerate}
		\item    $M$ is diffeomorphic to  one of the following four surfaces: a cylinder $S^1\times [0,1]$, a disk $D^2,$ a sphere $S^2$, or a torus $T^2$; 
		\item a function $f$ has no saddles, each isolated local extreme of $f$ is non-\textr{degenerate}, and the number of such local extremes is equal to $\chi(M)$.
		\item  if $M$ is diffeomorphic to $T^2$, then $f$ has  an even number of extremal circles. In particular, if $f$ is null-homotopic, then it has at least two extremal  circles.
	\end{enumerate}
\end{proposition}

These results are  direct generalizations of known results on the homotopy type of $\Sid$ for functions from $\F$ with {\it  only isolated singularities}; see   \cite[Theorem 3.7]{Maksymenko:OsakaJM:2011}, and  \cite[Theorem 1.3]{Maksymenko:AGAG:2006}.

\subsection{Structure of the paper}
The further text is organized in 11 sections.
In Section \ref{sec:Preliminaries},  we discuss  topological structure of functions from the class $\F$ near their critical points (\textsection \ref{subsec:hom-poly} and \ref{ssec:class-F}).
The foliation $\Delta_f$ on $M$ induced by functions from $\F$ is introduced in \textsection \ref{ssec:f-adapted}.

Section \ref{sec:SS} contains some generalities on  flows of vector fields on surfaces (\textsection\ref{ssec:generalit-flow}) and on shift maps along their trajectories (see \textsection\ref{ssec:shift-map-def}).
Diffeomorphisms of $\bbR^2$ that are shifts along trajectories of vector fields of the form $F_{\mu,n} = \mu y^n\frac{\partial}{\partial x}$ on $\bbR^2$, where $\mu:\bbR^2\to \bbR$ is a positive smooth function, $n\geq 0$,  are studied in Section \ref{sec:hor-lines}.

In Section \ref{sec:H-fileds}, we consider special vector fields on surfaces, called $H$-fields, naturally associated with a function from $\F$. In particular,  some relevant facts on Hamiltonian vector fields will be recalled in \textsection\ref{ssec:Ham-skew} and \ref{ssec:ham-function-F}.
In \textsection\ref{ssec:H-like}, we introduce the notion of an $H$-field for functions from $\F$ (see Proposition \ref{prop:vector-field-f}), which will be mainly used throughout the paper.
Proposition  \ref{prop:vector-field-f} is proved in Section \ref{sec:existence}.
 Section \ref{sec:Hom-type-DidF} is devoted to the study of the group of diffeomorphisms preserving flows of $H$-fields for functions from $\F$ and their homotopy properties.
The question of the existence of shift functions with respect to  flows of  $H$-fields for diffeomorphisms from $\Sid$ is discussed in Section \ref{sec:shift-diff}.

In Sections~\ref{sec:fibration-rho} and~\ref{sec:pi0G}, we obtain several auxiliary results on the homotopy properties of~$\mathcal{S}_{\mathrm{id}}(f)$ and its subgroups.
Finally, we prove Theorem~\ref{thm:main-theorem}  in Section~\ref{sec:proof-main}, and Proposition~\ref{prop:main} in Section~\ref{sec:proof-main-prop}.

\subsection{Strategy of the proof of Theorem \ref{thm:main-theorem}} 
The purpose of this paragraph is to outline the sketch of the proof of Theorem~\ref{thm:main-theorem}.

Let $M$ be a smooth, compact, connected, and \textr{orientable} surface, let $f$ be a function from $\F$
with the set of extremal circles $E_f$
 (see \textsection\ref{ssec:class-F}). Denote by $\mathcal{D}(M,E_f)$  the group of \textr{diffeomorphisms} of $M$ fixed on $E_f$.

{Roughly speaking, the homotopy type of $\Sid$ is determined by the homotopy properties of  two its subgroups
$\mathcal{G}(f,E_f)$ and $\mathcal{S}_{\id}(f,E_f)$, where $\mathcal{S}_{\id}(f,E_f)$ is a \textr{path} component of $\mathcal{S}(f,E_f) = \mathcal{S}(f)\cap \mathcal{D}(M,E_f)$ containing $\id_M$ and $\mathcal{G}(f,E_f) = \mathcal{D}(M,E_f)\cap \Sid$, see Eq.~\eqref{eq:Sid-inclusions}.}

Note that, in general, the group   $\mathcal{G}(f,E_f)$  is  not  connected, whereas  $\mathcal{S}_{\id}(f, E_f)$ coincides with the \textr{path} component  $\mathcal{G}_{\id}(f,E_f)$ of $\mathcal{G}(f,E_f)$ containing $\id_{M}$.
\textr{ If $E_f = \varnothing$, then all three groups $\mathcal{S}_{\id}(f)$, $\mathcal{G}(f,E_f)$ and $\mathcal{S}_{\id}(f,E_f)$ coincide; see Lemma \ref{lm:inclusions-stab}.}

For convenience, the proof may be divided into several steps. 
The goal for first two steps is to describe the homotopy type of the group $\mathcal{S}_{\id}(f,E_f)$. This can be done similarly to the proof of \cite[Theorem 1.3]{Maksymenko:AGAG:2006}.

{\it Step 1.}  
For $f\in \F$, we define a vector field $F$ (Proposition \ref{prop:vector-field-f}), called an $H$-like field associated with $f$, which encapsulates the necessary properties of the foliation $\Delta_f$ on $M$ induced by $f$ (see \textsection\ref{ssec:f-adapted}). Denote by $\bfF$ the flow of an $H$-like field of $f$.
{The homotopy type of the group $\mathcal{D}_{\id}(\bfF)$ of  diffeomorphisms of $M$ preserving trajectories of $\bfF$ and  isotopic to $\id_M$ via an isotopy which also preserves orbits of $\bfF$ is computed in  \cite[Theorem 3.5]{Maksymenko:OsakaJM:2011} (see Lemma \ref{lm:hom-type-Did}) --- $\mathcal{D}_{\id}(\bfF)$  is either contractible, or has the homotopy type of $S^1$ (see Section \ref{sec:Hom-type-DidF} for further discussion).}

{\it Step 2.}
{Then we  show (Proposition \ref{prop:shift-functions-exist-global}) that $\mathcal{S}_{\id}(f,E_f) = \mathcal{D}_{\id}(\bfF)$. Therefore, $\mathcal{S}_{\id}(f,E_f) $ is either contractible or has the homotopy type of a circle $S^1$; see Corollary \ref{cor:hom-type-SidE} for the precise statement.

Recall that if $E_f = \varnothing$, then $\mathcal{S}_{\mathrm{id}}(f) = \mathcal{S}_{\mathrm{id}}(f,\varnothing)$. In this case, Proposition \ref{prop:shift-functions-exist-global} states that $\mathcal{S}_{\mathrm{id}}(f) = \mathcal{D}_{\mathrm{id}}(\mathbf{F})$.}
 This  completes the proof of Theorem \ref{thm:main-theorem} for functions from $\F$ without extremal circles. The above arguments include, as a partial case, the proof of \cite[Theorem 1.3]{Maksymenko:AGAG:2006} for  functions on \textr{orientable} surfaces.

\textr{{\it From now on, we assume that the set of extremal circles $E_f$ of $f$ is not empty, i.e., $E_f = \{ C_1,\ldots, C_n\}$ for some $n\geq 1$.} }

{\it Step 3}. At this step, we establish some results on homotopy properties of $\Sid$.
 O.~Khohliyk and S.~Maksymenko proved  \cite[Theorem 8.2., Theorem 3.3.]{KhohliykMaksymenko:Indag:2020} (see Theorem \ref{thm:fibration-rho}) that there exists a locally trivial fibration $\rho_0:\Sid\to \mathcal{D}_{\id}(E_f)$ with the  fiber $\mathcal{G}(f,E_f)$, where $\mathcal{D}_{\id}(E_f)$ is a \textr{path} component of the identity map of the group $\mathcal{D}(E_f)$ of diffeomorphisms of an $1$-manifold $E_f$.
 
 Using the local triviality of $\rho_0$ and classical results of R.~Palais~\cite{Palais:MemAMS:1960} and J.~Milnor~\cite{Milnor:Trans:1959}, we prove (Lemma \ref{lm:CW}) that $\Sid$ has the homotopy type of a CW complex.
 
\textr{From the fact that $\mathcal{D}_{\id}(E_f)$ is homotopy equivalent to a torus $T^{|E_f|}$ (see Lemma \ref{lm:DidE}), where $|E_f|$ is the number of path components of $E_f$ (cardinality of $E_f$),} together with a long exact sequence of homotopy groups of $\rho_0$ and the homotopy type of $\mathcal{S}_{\id}(f,E_f)$, we  obtain that $\pi_q\mathcal{S}_{\id}(f) = 0$ for $q\geq 2$ and $\pi_1\Sid$ is the part of the following short exact sequence:
\begin{equation}
\xymatrix{
	1\ar[r] & \pi_1\Sid \ar[r]^{~~(\rho_0)_*} & \bbZ^{|E_f|}\ar[r] & \pi_0\mathcal{G}(f,E_f)\ar[r] &1,
} \tag{\ref{eq:non-triv-subseq}}
\end{equation}
see Lemma \ref{lm:weak-Sid}.
By Whitehead theorem, it follows that the homotopy type of $\Sid$ is completely determined by its fundamental group $\pi_1\Sid$.

{\it Step 4.} In order to describe $\pi_1\Sid$, one needs to study the group $\pi_0\mathcal{G}(f, E_f)$.
 In Section \ref{sec:pi0G},
we prove (Proposition \ref{prop:pi0-G}) that   $\pi_0\mathcal{G}(f,E_f)$ is a free abelian group of the rank $|E_f|$  or $|E_f|-1$. Thus, a short exact sequence \eqref{eq:non-triv-subseq} always splits, and consequently,
$\pi_1\Sid$ is either a trivial group, or it is isomorphic to  $\bbZ$. Therefore,
 by Whitehead theorem, $\Sid$
 is either  contractible or  homotopy equivalent to $S^1$.

\subsection*{Conventions} Throughout the paper we  work in the category of smooth \textr{($C^{\infty}$-differentiable)} manifolds. Unless stated otherwise, all objects naturally associated with a smooth surface $M$ -- such as functions, vector fields, their flows, etc. -- are assumed to be \textr{$C^{\infty}$-differentiable}; all subspaces of $C^{\infty}(M,M)$ and $C^{\infty}(M,P)$ are endowed with the subspace topologies. 
Moreover, the surface $M$ is always assumed to be compact, connected, and \textr{orientable}. For a topological group $G\in C^{\infty}(M,M)$, we denote by $G_{\id}$ the identity \textr{path} component of $G$.

\subsection*{Acknowledgment} 
This work was supported by a grant from the Simons Foundation (SFI-PD-Ukraine-00014586 B.G.F). 

The author would like to express his gratitude to Yevgen Polulyakh for pointing out an error in the first version of the proof of Proposition \ref{prop:vector-field-f}, and to Sergiy Maksymenko for valuable discussions.

The author is also grateful to the anonymous Referee for a careful reading of the manuscript and helpful suggestions which allowed to clarify the exposition.

\section{Generalities on class $\mathcal{F}$}\label{sec:Preliminaries}
Let $M$ be a smooth, connected,  compact and \textr{orientable} surface, and let $f$ be a function from $\F$ with the set of critical point $\Crit$. The set $\pi_0 \Crit$ of connected components  of $\Crit$
is the union of two sets $\CritC$ and $\CritP$, where $\CritC$ is the set of  isolated critical points and $\CritC$ is the set of critical circles of $f$.

By our definition, isolated singularities of $f\in \F$ are ``modeled'' by homogeneous polynomials without multiple factors. Therefore, in the next paragraph, we discuss foliations on $\bbR^2$ by the level-sets of such polynomials.
\subsection{Homogeneous polynomials}\label{subsec:hom-poly}
\begin{figure}[h]
	\centering
	\includegraphics[width=0.6\textwidth]{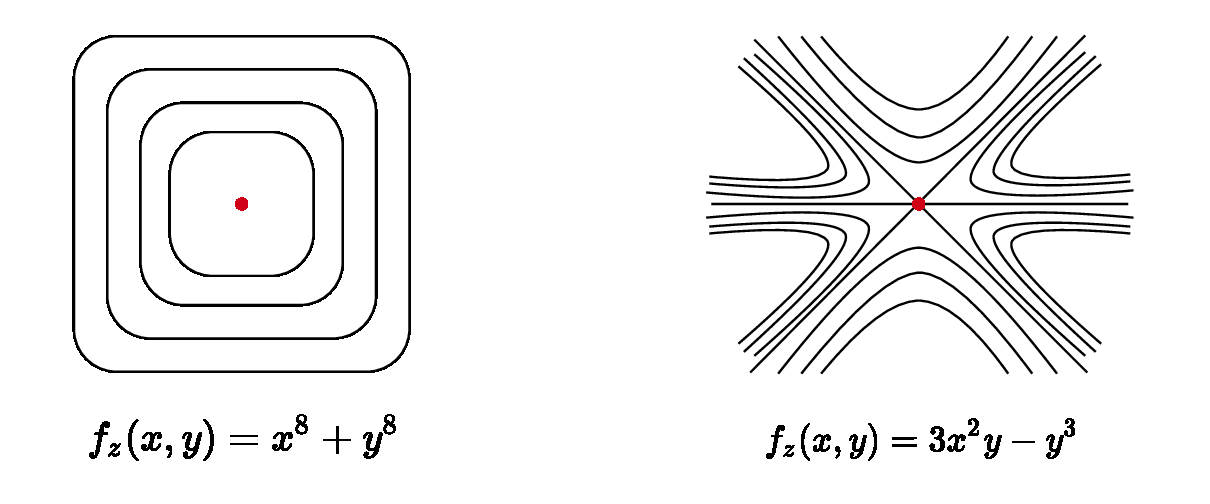}
	\caption{Foliations by homogeneous polynomials: a \textr{degenerate} extreme (left) and a  $3$-saddle (right)}
	\label{fig:crit-points}
\end{figure}

 Let $f:\bbR^2\to \bbR$ be a real homogeneous polynomial. It is well known that $f$ factors over $\bbR$ into a finite product of linear $L_i = a_ix+b_iy$ and irreducible over $\bbR$ quadratic factors $Q_j(x,y) = c_jx^2+2d_jxy+e_jy^2$, i.e.,
\begin{equation}\label{eq:homogeneous}
f(x,y) = \prod_{i = 1}^p L_i(x,y)\cdot \prod_{j=1}^q Q_j(x,y),
\end{equation}
for some $p,q\geq 1$.
The origin $0\in \bbR^2$ is the only critical point of \eqref{eq:homogeneous} if and only if $\deg f\geq 2$ and $f$ has no multiple factors.
For such polynomials, the origin $0\in \bbR^2$ will be called
\begin{itemize}
	\item  an {\it extreme}, if $f = Q_1$ (non-\textr{degenerate}), or $f = Q_1Q_2\ldots Q_q$ (\textr{degenerate}),
	\item a {\it saddle}, if $f = L_1 Q_1Q_2\ldots Q_q$ (quasi-saddle), or $f = L_1L_2$ (non-\textr{degenerate} or $2$-saddle), or if $\deg f = p+2q\geq 3$ for $p\geq 2$ (generalized saddle or $p$-saddle).
\end{itemize}
Examples of homogeneous polynomials and their level sets are shown in Fig.~\ref{fig:crit-points}.

Let us note that isolated singularities of functions from $\F$ are ``topologically generic'',
 since    by result of P.~T.~Church and J.~G.~Timourian \cite{ChurchTimourian}, and  A.~Prishlyak \cite{Prishlyak:2002},  the topological structure of level sets of any smooth function on surface near an isolated critical point
can be realized by the level sets of some homogeneous polynomial without multiple factors; see also the discussion in \cite{Maksymenko:DefFuncI:2014}.

\subsection{Critical circles of $f$}\label{ssec:class-F}
Let $C\in \CritC$ be a critical circle of $f$ and $z\in C$. By (2b) of Definition \ref{def:class-F}, the germ at $z$ of $f$  is smoothly equivalent to a germ at $0\in \mathbb{R}^2$ of $f_C(x,y) = \pm y^{n_C}$  for  $n_C\geq 2$ depending on $C.$ Note that we allow the number $n_C$ to be either odd or even. 

\begin{figure}[h]
	\centering
	 \includegraphics[width=0.45\textwidth]{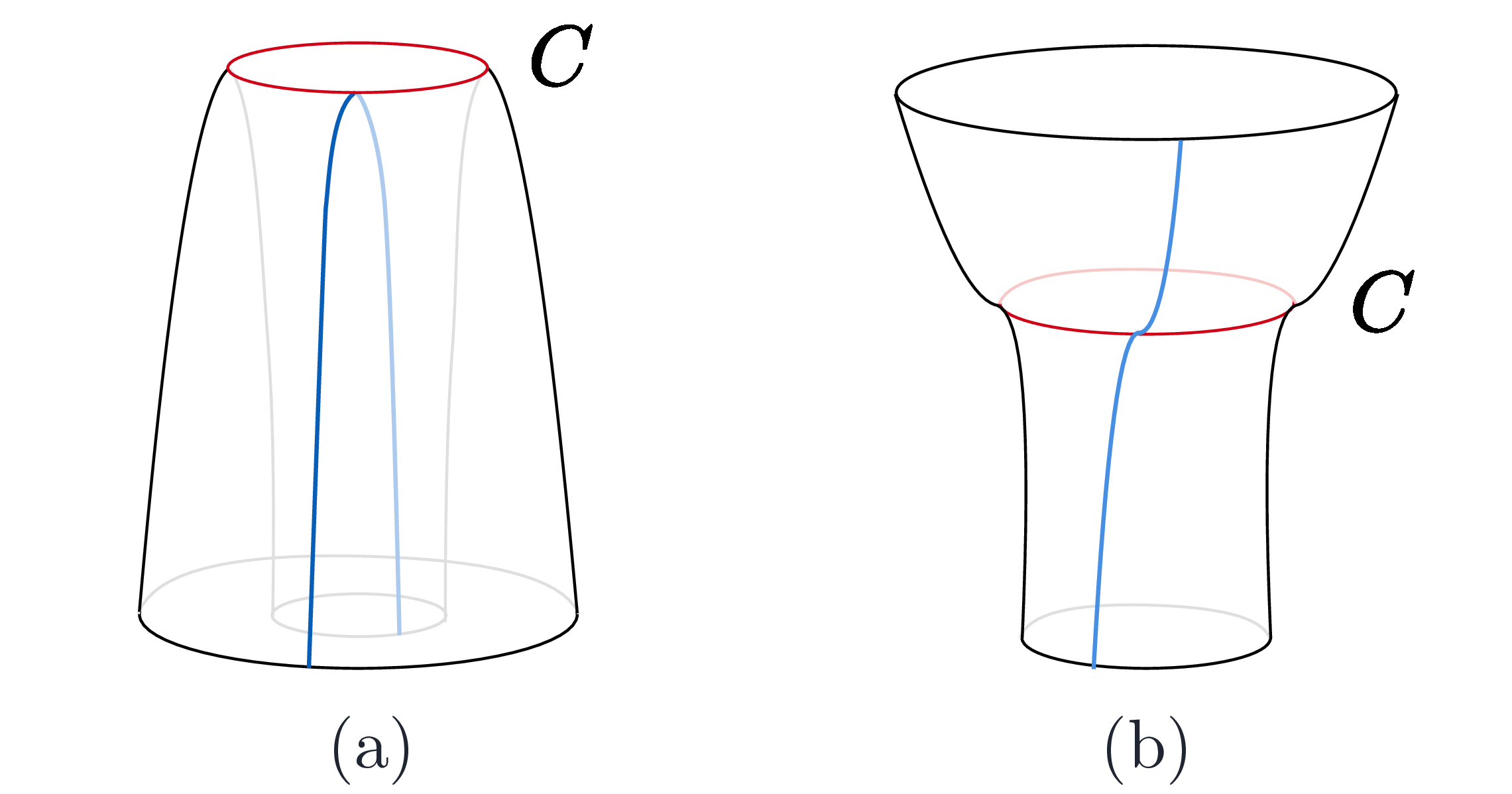}
	 \caption{Neighborhood of a critical circle $C$ on an \textr{orientable} surface: (a) $n_C$ is even, (b) $n_C$ is odd}
	 	\label{fig:crit-circles}
\end{figure}

Let $\mathbf{G}:M\times \bbR\to M$ be the flow of the gradient vector field $\mathrm{grad}(f)$ of $f$ with respect to some Riemannian metric on $M$. 
Let $W$ be an open and  connected  neighborhood of $C$ \textr{which is a union of} connected component of level-sets of $f$ and such that $W\setminus C$  contains no critical points of $f$.  Since $M$ is an \textr{orientable} surface, it follows that $W$ is a cylinder.
If $n_C$ is even, then $C$ is \textr{an} {\it extremal} critical circle for $f$, so either $\lim\limits_{t\to \infty}\mathbf{G}_t(p)\in C$ ($C$ is {\it maximal}) or $\lim\limits_{t\to- \infty}\mathbf{G}_t(p)\in C$ ($C$ is {\it minimal}) for $p\in W$. If $n_C$ is odd, then $C$ is {\it non-extremal},  see Fig. \ref{fig:crit-circles}. 

It is easy to see that the class $\F$ contains the class of Morse-Bott function and the class of  Morse functions on $M$ that are locally constant on $\partial M$. So the class $\F$ is ``massive'' and consists of functions with ``topologically generic'' singularities.

\subsection{Saturated neighborhoods}\label{ssec:f-adapted}
A function $f$ from $\F$ induces a  foliation $\Delta_f$ with singularities on $M$.
A leaf $K$ of $\Delta_f$ is either an isolated critical point or a connected component of the complement $f^{-1}(c)\setminus \CritP$, $c\in P.$  Thus, the only singular  leaves of $\Delta_f$ are the isolated critical points of $f$, while a critical circle  $C$ of  $f$ is a regular leaf of $\Delta_f$.

A subset $X\subset M$ is called \textr{{\it saturated}} if $X$ consists of leaves of $\Delta_f$.
 An open and connected neighborhood $U$ of a leaf $\gamma\in \Delta_{f}$ which consists of regular leaves of $\Delta_f$ will be called a \textr{{\it saturated neighborhood}} of $\gamma.$

\section{Smooth shifts along trajectories of flows}\label{sec:SS}
\subsection{Generalities}\label{ssec:generalit-flow}
Let $M$ be a smooth and compact surface, and let $F$ be a vector field on $M$ tangent to $\partial M$ with the flow $\mathbf{F}:M\times\bbR\to M$. Denote by $\Sigma_F$  \textr{the} set zeros of $F$. For a periodic (w.r.t the flow $\bfF$) point $x\in M$, we denote by $\mathrm{per}_{\bfF}(x)$ its prime period. Note that if all regular points of $\bfF$ are periodic, then
  $\mathrm{per}_{\bfF}:M\setminus \Sigma_F\to \bbR$ is a smooth function which is constant along the orbits of $\bfF$.
  
A flow $\bfF$ defines a natural  foliation $\Delta_{\bfF}$ with singularities on $M$, whose leaves are trajectories of $\bfF.$  Denote by $\mathcal{D}(\bfF)$  the group of diffeomorphisms of $M$, which preserve the leaves of $\Delta_{\bfF}$, i.e, for $h\in \mathcal{D}(\bfF)$ we have $h(\gamma)\subset \gamma$ for all $\gamma\in \Delta_{\bfF}$. Let also $\mathcal{D}^+(\bfF)$ be a subgroup of $\mathcal{D}(\bfF)$ consisting of   diffeomorphisms of $M$  preserving  orientation of regular trajectories of $F$, and $\mathcal{D}_{\id}(\bfF)$ be the \textr{path} component of $\mathcal{D}(\bfF)$ containing  $\id_M$.

Let $\mu:M\to \bbR$ be a smooth function and $\mathbf{G}:M\times \bbR\to \bbR$ be the flow of the vector field $G = \mu F.$ It is known \cite[Lemma 2.1]{Maks:reparam-sh-map} that there exist a smooth function $\sigma:M\to \bbR$ such that 
\begin{equation}\label{eq:repar-flow}
\mathbf{G}(x,t) = \mathbf{F}(x,\sigma(x,t)), \qquad \sigma(x,t)  = \int_0^t\mu(\mathbf{G}(x,\tau))d\tau.
\end{equation}
Hence, the flow $\mathbf{G}$ is a {\it reparametrization} of the flow $\bfF$.  
Note that if $\mu\neq 0$ on $M$, then $\bfF$ and $\mathbf{G}$ define the same foliation on $M$.

Let $V$ be a subset of $M$ such that $V$ contains no zeros of $F$. {Two vector fields $G$ and $F$ are said to be {\it codirectional} on $V$, if there exists a smooth function $\mu:V\to \bbR$ such that  $G = \mu F$ and $\mu|_V>0$.}

\subsection{Shifts along trajectories of flows}\label{ssec:shift-map-def}
Let $U$ be an open subset of $M$. For a smooth function $\alpha:U\to \bbR$, define the following smooth map
\begin{equation}
\bfF_{\alpha}:U\to M,\qquad\bfF_{\alpha}(x) = \bfF(x,\alpha(x)), \qquad x\in U.
\end{equation}
Note that for each $\alpha\in C^{\infty}(U, \bbR)$, the map $\bfF_{\alpha}$ preserves the orbits of $\bfF$ on $U$, 
i.e., 
$\bfF_{\alpha}(\gamma\cap U)\subset \gamma$, where $\gamma\in \Delta_{\bfF}$. The map
$\bfF_{\alpha}$
is called a {\it shift along $\bfF$}.
Roughly speaking,  under the action of $\bfF_{\alpha}$, the point $x$  shifts along its trajectory of  $\bfF$  by its ``own time''  $\alpha(x)$, which depends smoothly  on $x$.

Let $h:U\to M$ be a smooth map which preserves trajectories of $\bfF$ on $U,$ i.e., $h(\gamma\cap U)\subset \gamma$ for all  $\gamma\in \Delta_{\bfF}$. 
We say that $h$ is a {\it shift along $\bfF$} on $U$ if  $h = \bfF_{\alpha}$ for some smooth function $\alpha:U\to \bbR$. Such smooth function $\alpha$ is called {\it a shift function} for $h$ on $U$ (with respect to $\bfF$).
Some basic properties of such maps are stated in the following lemma.

\begin{lemma}[Lemma 2, Proposition 3, Theorem 19 in \cite{Maksymenko:TA:2003}; Lemma 2.3 \cite{Maksymenko:OsakaJM:2011} or Lemma 6.1 \cite{Maksymenko:DefFuncI:2014}]\label{lm:shift-main}
	Let $U$ be an open subset of $M$.
	The following statements hold.
	\begin{enumerate}\label{lm:prop-sh-func}
		\item Let $h:U\to M$ be a smooth map such that $h(\gamma\cap U)\subset \gamma$ for each orbit $\gamma$ of $\bfF$, $z\in U$ be non-zero point of $F$, and $a\in \bbR$ be a number such that $h(z) = \bfF(z,a)$. Then there exists an open neighborhood $V$ of $z$ and a unique smooth function $\alpha:V\to \bbR$ with $\alpha(z) = a$ and such that  $h(x) = \bfF(x,\alpha(x))$ for all $x\in V.$ 
		\item Assume that $U$ does not contain zeros of $F$. Let 
		  $\alpha,\beta:U\to \bbR$ be smooth functions such that $\bfF_{\alpha} = \bfF_{\beta}$ on $U$. If $\alpha(z) = \beta(z)$ for some $z\in U$, then $\alpha = \beta$ on $U$. In particular, if $U$ contains a non-closed orbit of $\bfF$, then $\alpha = \beta$ on $U$.
		\item Consider the following subset 
		$$
		\Gamma_{U,F} = \{\alpha\in C^{\infty}(U,\bbR),|\, 1+F\alpha\neq  0  \}
		$$ 
		 of $C^{\infty}(U,\bbR)$.  Then $h = \bfF_{\alpha}:U\to \bfF_{\alpha}(U)$ is a diffeomorphism iff $\alpha\in \Gamma_{U,F}$.
		\item for any smooth functions $\beta:U\to \bbR$ and $\alpha:\bfF_{\beta}(U)\to \bbR$, we have
		$\bfF_{\alpha}\circ \bfF_{\beta} = \bfF_{\beta +\alpha\circ \bfF_{\beta}}$.
	\end{enumerate}
\end{lemma}

\subsection{Shift map}
Let $V$ be an open subset of $M$.
\textr{The following map}
$$
\phi_{V,\bfF}:C^{\infty}(V,\bbR)\to C^{\infty}(V,M), \qquad \phi_{V,\bfF}(\alpha) = \bfF_{\alpha}.
$$
is called a {\it  shift map} (along trajectories of $\bfF$) on $V$. The image $\mathrm{Im}(\phi_{V,\bfF})$  will be also denoted by $\mathrm{Sh}(V,F).$
The set $\ker\phi_{V,\bfF} = \phi^{-1}_{V,\bfF}(\id_M)$ is called the {\it kernel of} $\phi_{V,\bfF}.$ It is obvious that $0\in\ker\phi_{V,\bfF}$.

\begin{lemma}[Theorem 12, Proposition 14 in \cite{Maksymenko:TA:2003}, Theorem 1.1 \cite{Maks:reparam-sh-map}]	\label{lm:ker-shift}
	The following \textr{statements  hold.}
\begin{enumerate}
	\item \textr{Let $\alpha, \beta:V\to \bbR$ be smooth functions.} Then  	$\phi_{V,\bfF}(\alpha) = \phi_{V,\bfF}(\beta)$ iff $\alpha-\beta\in \ker\phi_{V,\bfF}.$
	\item A map $\phi_{V,\bfF}$ is locally injective iff $\Sigma_F\cap V$ is nowhere dense in $V$.
	\item  Assume that $V$ is connected and  $\Sigma_F\cap V$ is nowhere dense in $V$.  Then 
	\begin{enumerate}[label={\rm(\roman*)}]
		\item either $\ker\phi_{V,\bfF} = \{0\}$ and thus, $\phi_{V,\bfF}$ is injective. This case holds if $V$ contains  either a non-periodic point $F$, or a fixed point $z\in \Sigma_F\cap V$ such that the tangent flow $T_z\bfF_{t}$ on $T_zM$ is the identity,
		\item or $\ker\phi_{V,\bfF} = \{n\theta\}_{n\in \bbZ}$ for some smooth  function $\theta:V\to(0,\infty)$. In this case 
		 $\phi_{V,\bfF}$ yields a bijection between $C^{\infty}(M,\bbR)/\ker(\phi_{V,\bfF})$  and $\mathrm{Sh}(V, F)$. Therefore, for every $\alpha\in C^{\infty}(V,\bbR)$, we have
		$$
		\phi_{V,\bfF}^{-1}\circ \phi_{V,\bfF}(\alpha) = \alpha+\ker\phi_{V,\bfF} = \{\alpha+n\theta\}_{n\in \bbZ}.
		$$
		In particular, every non-zero point $z$ of $F$ on $V$ is periodic for some period $\mathrm{per}_{\bfF}(z)$, and
		$\theta(z)  = n\mathrm{per}_{\bfF}(z)$ for some $n\in \bbN$. Hence,  $\theta$ is constant along trajectories of $F$. Such  function $\theta$ is  called a {\it period function} for $\phi_{V,\bfF}.$ 
	\end{enumerate}
	\item Let $\mu:M\to \bbR$ be a smooth function. Then $\mathrm{Sh}(V,\mu F)\subset \mathrm{Sh}(V, F)$ and, in particular, if $\mu\neq 0$ on $M$, we have $\mathrm{Sh}(V,\mu F) = \mathrm{Sh}(V, F)$. 
\end{enumerate}
\end{lemma}
\noindent When $V = M$, we shall use the simplified notation $\phi_{\bfF}$ for $\phi_{M, \bfF},$ $\mathrm{Sh}(F)$ for $\mathrm{Sh}(M,F)$, and so on.
We also denote by $\mathcal{D}^{\mathrm{sh}}(F)$ the group of diffeomorphisms from $\mathcal{D}^+(\bfF)$ which are shifts along trajectories of  $F$, i.e., $\mathcal{D}^{\mathrm{sh}}(F) = \mathrm{Sh}(F)\cap \mathcal{D}^+(\bfF).$

\section{Singular foliations on the plane on horizontal lines}\label{sec:hor-lines}
In this paragraph, we recall some  properties of vector fields on $\bbR^2$ of the form $F_{\mu,n} = \mu (x,y) y^n\frac{\partial}{\partial x}$, where  $\mu:\bbR^2\to \bbR$ is a smooth positive function and $n\in \bbN\cup \{0\}$. Let $\bfF^{\mu,n}:\bbR^2\times \bbR\to \bbR^2$ be the flow of $F_{\mu,n}$. For $\mu = 1$, there is a simple explicit formula for $\bfF^{1,n}$:
\begin{equation}\label{eq:flow-yn}
\bfF^{1,n}(x,y,t) = (x+y^nt,y).
\end{equation}
{Note that for each $n\geq 1$,  the flow  $\bfF^{\mu,n}$ is a reparameterization of the flow  $\bfF^{1,n}$, see  Eq. \eqref{eq:repar-flow}.}

\begin{figure}[h]
	\centering
	\includegraphics[width=0.75\textwidth]{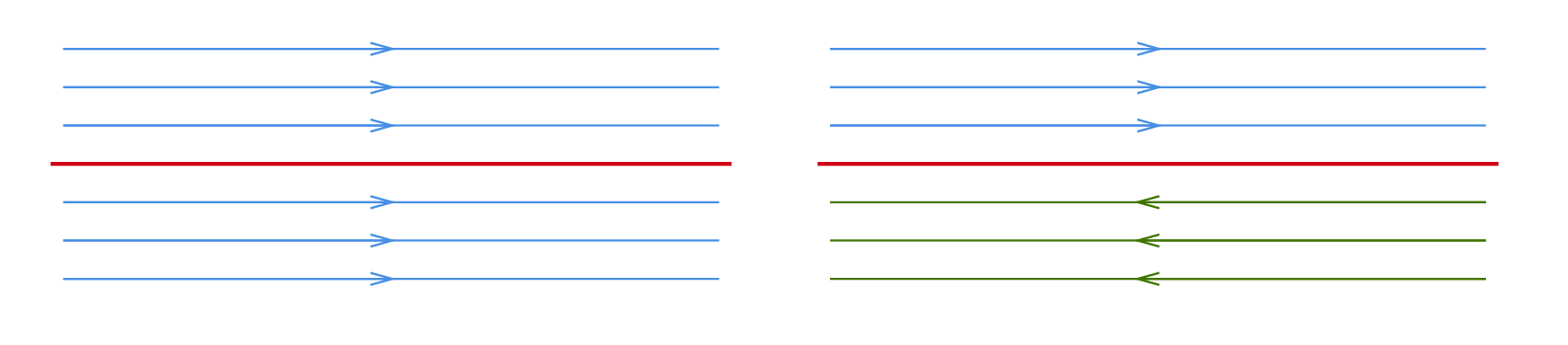}
	\caption{Trajectories of $F_{\mu,n}$, $n\geq 1$ when $n$ is even (left) and $n$ is odd (right)}
	\label{fig:trajectories}
\end{figure}

It is easy to see that $\bfF^{\mu,0}$ is non-singular, and the leaves of $\Delta_{\bfF^{\mu,0}}$ are horizontal lines.
If $n\geq 1$, the flow $\mathbf{F}^{\mu,n}$ has fixed points. 
The leaves of $\Delta_{\bfF^{\mu,n}}$ are the lines $\{y = a\}$, $a\in\bbR\setminus\{0\}$  (regular leaves) and the points $(x,0)$, $x\in\bbR$ (singular leaves).
Note that regular  trajectories of $\bfF^{\mu,n}$ may ``change'' their orientation when passing through the set of zeros of $F_{\mu,n}$ in the traversal direction; see
 Fig. \ref{fig:trajectories}.

{It is easy to prove that the following inclusions hold:
\begin{equation}\label{eq:D+F}
	\mathcal{D}_{\id}(\bfF^{\mu,n})\subset\mathcal{D}^{+}(\bfF^{\mu, n}) = \mathcal{D}(\bfF^{\mu,n})\cap \mathcal{D}^+(\bbR^2),
\end{equation}
where $\mathcal{D}^+(\bbR^2)$ is the group of orientation-preserving diffeomorphisms of $\bbR^2$.}

{Let $(u,v)$ be any smooth coordinate system on $\bbR^2$. With respect to $(u,v)$ denote by $F_{(u,v)}$ a vector field on $\bbR^2$ of the form 
$$
F_{(u,v)} = v\frac{\partial}{\partial u}.
$$
The following lemma is technical and it is needed for  proofs in Section \ref{sec:existence}.
\begin{lemma}\label{lm:technical-lema}
	Let $h:\mathbb{R}^2\to \bbR^2$ be a diffeomorphism given by $h(x,y) = (u(x,y), v(y))$ with $u_x >0$, $v'>0$ and $v(0) \equiv 0$. Denote by  $G = h^*(F_{(u,v)})$  a pull-back of $F_{(u,v)}$ along $h$. Then 
	\begin{enumerate}
		\item $G = \tau F_{(x,y)}$ for some positive smooth function $\tau:\bbR^2\to \bbR$,
		\item for any smooth  function $\rho:\bbR^2\to \bbR$ with $0\leq \rho\leq 1$ we have 
		$$
		\rho F_{(x,y)}+(1-\rho) G = \eta F_{(x,y)}
		$$
		 for some smooth and positive function $\eta:\mathbb{R}^2\to \bbR$.
	\end{enumerate}
\end{lemma}

\begin{proof}
(1)	Note that, since $v(0) \equiv 0$, it follows from Hadamard lemma that 
	\begin{equation}\label{eq:v-Hadamard}
		v(y) = ys(y)
	\end{equation}
	for each $y\in \mathbb{R}$, where $s:\mathbb{R}\to \mathbb{R}$ is a smooth and positive function. 

	Components of the vector field $G = h^*(F_{(u,v)}) =  \xi_1(x,y)\frac{\partial}{\partial x}+\xi_2(x,y)\frac{\partial}{\partial y}$  can be computed from the following diagram
	 $$
	 \xymatrix{
	T\bbR^2 \ar[r]^{Th} & T\bbR^2\\
	\bbR^2\ar[r]^h \ar[u]^{G} 	 &	\bbR^2\ar[u]_{F_{(u,v)}}
}
	 $$
	Direct calculation with 
	$$
	Th = \begin{bmatrix}
		u'_x& u'_y\\
		0& v'
	\end{bmatrix}
\qquad\text{and}\qquad F_{(u,v)}(h(x,y)) = \begin{bmatrix}
	v(y)\\0
\end{bmatrix}
	$$
yields $\xi_1(x,y) = \frac{v(y)}{u'_x(x,y)}$ and $\xi_2(x,y) = 0$, $(x,y)\in\bbR^2$, thus $G = \frac{v(y)}{u'_x(x,y)}\frac{\partial}{\partial x}.$ Applying \eqref{eq:v-Hadamard} to the formula for $G$, we obtain
	$$
	G =  \frac{v(y)}{u'_x(x,y)}\frac{\partial}{\partial x} = \frac{ys(y)}{u'_x(x,y)}\frac{\partial}{\partial x} = \frac{s(y)}{u'_x(x,y)}\Big(y\frac{\partial}{\partial x}\Big) = \tau F_{(x,y)},
	$$
	where $\tau:\bbR^2\to \bbR$, $\tau(x,y) = \frac{s(y)}{u'_x(x,y)}$ is a  smooth function. It is also positive since $u_x'>0$ and $s>0$.
	
	(2) is a direct consequence of (1). It is easy to see that the function $\eta$ has the form  $\eta = (\rho + (1-\rho)\tau)$.
\end{proof}
}

\subsection{Shifts along $F_{\mu,n}$}
In this paragraph, we  study the relations between the groups $\mathcal{D}^{+}(\bfF^{\mu, n})$ and $\mathcal{D}^{\mathrm{sh}}(F_{\mu,n})$; see Section \ref{sec:SS} for the definitions.

In the simplest case, when $n = 0$, a vector field $F_{1,0} = \frac{\partial}{\partial x}$ has no zeros in $\bbR^2$ and  one can show using (1) and (2) of Lemma \ref{lm:shift-main} that for each $h\in\mathcal{D}^+(\mathbf{F}^{1,0})$ there exists a unique smooth function $\beta:\bbR^2\to \bbR$ such that $h = \mathbf{F}^{1,0}_{\beta}.$  
Since $\bfF^{\mu,0}$ and $\bfF^{1,0}$ define the same foliation on $\bbR^2$, it follows that $\mathcal{D}^+(\bfF^{\mu,0}) = \mathcal{D}^+(\bfF^{1,0})$.
By (4) of Lemma \ref{lm:ker-shift}, we finally get
$$
\mathcal{D}^{\mathrm{sh}}(F_{\mu,0}) = \mathcal{D}^{\mathrm{sh}}(F_{1,0}) = \mathcal{D}^+(\bfF^{1,0}) = \mathcal{D}^+(\bfF^{\mu,0}).
$$
If $n\geq 1$, then $\bfF^{\mu,n}$ is singular, and the situation is harder.

\begin{lemma}[cf. p.24 case (3) in \cite{Maksymenko:OsakaJM:2011}]\label{lm:example}
	Let $\mu:\bbR^2\to \bbR$ be a positive smooth function on $\bbR^2$ and $F_{\mu,n} = \mu y^n\frac{\partial }{\partial x}$ be a vector field on $\bbR^2$, $n\geq 1$.
	Then 
	
	{\rm (1)}
	the following inclusions hold:
	$$
	\xymatrix{
	\mathcal{D}^{\mathrm{sh}}(F_{\mu,n})	\ar@{=}[d] \ar@{^{(}->}[r]& \mathcal{D}^{\mathrm{sh}}(F_{\mu,n-1}) \ar@{=}[d] \ar@{^{(}->}[r]& \ldots \ar@{^{(}->}[r]& \mathcal{D}^{\mathrm{sh}}(F_{\mu,2}) \ar@{=}[d] \ar@{^{(}->}[r]& \mathcal{D}^{\mathrm{sh}}(F_{\mu,1}) \ar@{=}[d]  &  \\
	\mathcal{D}^{\mathrm{sh}}(F_{1,n})	\ar@{^{(}->}[d] \ar@{^{(}->}[r]& \mathcal{D}^{\mathrm{sh}}(F_{1,n-1}) \ar@{^{(}->}[d] \ar@{^{(}->}[r]& \ldots \ar@{^{(}->}[r]& \mathcal{D}^{\mathrm{sh}}(F_{1,2}) \ar@{^{(}->}[d]\ar@{^{(}->}[r]& \mathcal{D}^{\mathrm{sh}}(F_{1,1}) \ar@{^{(}->}[d] & \\
	\mathcal{D}^+(\bfF^{1,n}) \ar@{=}[r] \ar@{=}[d]& 	\mathcal{D}^+(\bfF^{1,n-1}) \ar@{=}[r] \ar@{=}[d]&  \ldots \ar@{=}[r] &	\mathcal{D}^+(\bfF^{1,2}) \ar@{=}[r] \ar@{=}[d]& 	\mathcal{D}^+(\bfF^{1,1}) \ar@{=}[d]& 	\\
	\mathcal{D}^+(\bfF^{\mu,n}) \ar@{=}[r]& 	\mathcal{D}^+(\bfF^{\mu,n-1}) \ar@{=}[r]&  \ldots \ar@{=}[r]&	\mathcal{D}^+(\bfF^{\mu,2}) \ar@{=}[r]& 	\mathcal{D}^+(\bfF^{\mu,1}) &
}
	$$
	
	{\rm (2)} $\mathcal{D}^{\mathrm{sh}}(F_{1,k}) = \mathcal{D}^+(\bfF^{1,k})$ iff $k = 1$.
\end{lemma}
Therefore, not any diffeomorphism from $\mathcal{D}^+(\bfF^{\mu,n})$ is a shift along trajectories of $\bfF^{\mu, n}$ for $n\geq 2$, but any $h\in \mathcal{D}^+(\bfF^{\mu,n})$ always has a shift function with respect to the flow $\bfF^{\mu,1}$.

\begin{proof}	
	(1)  Note that for any $\mu:\bbR^2\to \bbR$ and $1\leq k\leq n$, vector fields $F_{\mu,k}$ define the same foliation on $\bbR^2$, hence, $\mathcal{D}(\bfF^{\mu,k}) = \mathcal{D}(\bfF^{1,1})$. It follows from formula \eqref{eq:D+F} that $\mathcal{D}^+(\bfF^{\mu,k}) = \mathcal{D}^+(\bfF^{1,1})$.
	All other inclusions in the diagram above trivially follow from  Eq. \eqref{eq:repar-flow} and (4) of Lemma  \ref{lm:ker-shift}. We left details to the reader.

	(2) 
	It remains to prove that $\mathcal{D}^+(\bfF^{1,k})\subset \mathcal{D}^{\mathrm{sh}}(F_{1,k})$ only if $k = 1$. 
	Let  $h$ be a diffeomorphism  from $\mathcal{D}^{+}(\bfF^{1,k}),$ i.e.,    $h$ preserves  trajectories of $\bfF^{1,k}$ and their orientations. 
	Since all trajectories of $\bfF^{1,k}$ on $\bbR^2\setminus \{y = 0\}$ are non-closed and non-singular, then by (1) of Lemma \ref{lm:prop-sh-func}, a diffeomorphism $h$ has a shift function $\alpha:U\to \bbR$ on each open and connected  $U\subset \bbR^2\setminus \{y= 0\}$. By (2) of Lemma \ref{lm:prop-sh-func},  this function can be extended to some smooth function $\alpha:\bbR^2\setminus\{y = 0\}\to \bbR$ such that $h = \bfF^{1,k}_{\alpha}$ on $\bbR^2\setminus \{y= 0\}.$  

	Now we need to study the question of extending a shift function $\alpha$ to the whole plane $\bbR^2$.
	Using the  formula \eqref{eq:flow-yn} for the flow $\bfF^{1,k}$, we obtain
	\begin{equation}\label{eq:shift-alpha}
	h(x,y) = (h_1(x,y), h_2(x,y)) =\mathbf{F}^{1,k}(x,y,\alpha(x,y)) = (x+y^k\alpha(x,y), y)
	\end{equation}
	on $\bbR^2\setminus \{y=0\}.$
	Expressing $\alpha$ from Eq. \eqref{eq:shift-alpha}, we get
	\begin{equation}\label{eq:h_1}
		\alpha(x,y) = \frac{h_1(x,y)-x}{y^k}
	\end{equation}
	for $(x,y)\in \bbR^2\setminus \{y=0\}.$ A function $h_1(x,y)-x$ is  $C^{\infty}$-smooth and {$h_1(x,0)\equiv x$} for all $x\in \bbR$, since for any $x\in\bbR$ the point $(x,0)$ is a zero  of  $F_{1,k}$. By {Hadamard lemma, there exists a smooth function $\gamma:\bbR^2\to \bbR$ such that $h_1(x,y)-x = y\gamma(x,y)$}.
	Then Eq.~\eqref{eq:h_1} takes the following form
	$$
	\alpha(x,y) = \frac{y\gamma(x,y)}{y^k} = \frac{\gamma(x,y)}{y^{k-1}}
	$$
	for $(x,y)\in \bbR^2\setminus \{y = 0\}.$ 
	
	If $k = 1$, then $\alpha = \gamma$ on $\bbR^2\setminus \{y=0\}$, and since $\gamma$ is smooth on $\bbR^2$, it follows that $\alpha$ can be extended to a smooth function on $\bbR^2$. In particular, by Hadamard lemma, this function $\alpha:\bbR^2\to \bbR$ is given by
	 $$
	\alpha(x,y) =  \gamma (x,y) = \int\limits_0^1\frac{\partial h_1(x, ty)}{\partial y}dt.
	$$
	Hence,  we get $\mathcal{D}^+(\bfF^{1,1})\subset \mathcal{D}^{\mathrm{sh}}(F_{1,1})$.
\end{proof}

\section{$H$-fields for functions from $\F$}\label{sec:H-fileds}
Let $M$ be a smooth and {\it \textr{orientable}} surface.
It is well-known that $M$ admits a symplectic form $\omega:TM\times TM\to \bbR$, i.e., a non-\textr{degenerate} skew-symmetric $2$-form.
\subsection{Hamiltonian vector fields}\label{ssec:Ham-skew}
 From the fact that $\omega$ is non-\textr{degenerate}, \textr{it} follows that the map $\omega^{\flat}:TM\to T^*M$ given by $\omega^{\flat}_p(X) = \omega_p(X,\cdot)$,  $p\in M$ is a bundle isomorphism.  Let $f\in C^{\infty}_{\partial}(M,P)$ be a function and \textr{$Tf:TM\to TP$ be its tangent map.}
Since $P$ is either $\bbR$ or $S^1$, it is well-known that  $TP$ is a globally trivializable bundle; denote by $\zeta:TP\to P\times \bbR$  its trivialization isomorphism. Therefore, the map $Tf$  induces  \textr{a map $Df:TM\to \bbR$}, defined as the composition 
$$
\xymatrix{
Df:TM\ar[r]^-{Tf} & TP \ar[r]^-{\zeta}& P\times \bbR \ar[r]^-{p_2}& \bbR,
}
$$ 
where $p_2$ is the projection onto a second factor. \textr{Note that the map $Df$ is a usual differential of $f$ as a differential $1$-form.}
 A vector field $X_f = (\omega^{\flat})^{-1}(Df)$ on $M$ is called {\it a Hamiltonian vector field of $f$}. It is known that $X_f$ satisfies the following properties:
\begin{itemize}
	\item zeros of $X_f$ correspond to the critical points of $f$,
	\item $f$ is constant along the trajectories of $X_f$, i.e., $X_f(f) = 0.$ In other words, $X_f$ is tangent to the level sets of $f$ (and to $\partial M$).
\end{itemize}

\subsection{Hamiltonian vector fields for functions from $\F$}\label{ssec:ham-function-F}
For a function  $f\in\F$, we  have simple formulas for $X_f$ in some neighborhoods of critical points of $f$. In the chart $(U, (x,y))$ near $z\in\Crit$
from (2) of  Definition \ref{def:class-F}, a vector field 
$X_f$ on $U$ has the form
\begin{align}
	X_f &= \frac{1}{g(x,y)}  (f_z)'_y\frac{\partial}{\partial x} - \frac{1}{g(x,y)} (f_z)'_x\frac{\partial}{\partial y}, &(\text{if $z$ is isolated})\label{eq:ham-vf-near-iso}   \\
	 X_f &= \pm\frac{n_Cy^{n_C-1}}{{g}(x,y)}\frac{\partial}{\partial x},& (\text{if $z$ is not isolated, $z\in C$})\label{eq:Ham-near-C}
\end{align}
where $C$ is a critical circle of $f$,
$f_z$ and $n_C$ are such in (2) of  Definition \ref{def:class-F}, and $g:U\to \bbR$ is a positive smooth function such that $\omega = g(x,y)dx\wedge dy$ on $U$. Note that vector fields as \eqref{eq:Ham-near-C} were studied in \textsection\ref{sec:hor-lines}.

\subsection{$H$-fields of  functions from $\F$} \label{ssec:H-like}
The following proposition holds true.
\begin{proposition}\label{prop:vector-field-f}
Let $M$ be a smooth,  compact, connected, and \textr{orientable} surface, and let  $f$ be a function from $\F$ with the set of critical points $\Crit$. There exists a vector field $F$ on $M$ which satisfies:
\begin{enumerate}
	\item $Ff = 0$, i.e., $f$ is constant along trajectories of $F$;
	\item $F$ has no zeros in $M\setminus \Crit$;
	\item let $C$ be a connected component of $\Crit$, and 
	 $z\in C$ be a critical point of $f$. There exists a coordinate chart $(U_z,x,y)$ near $z$ such that $z = (0,0)\in U_z$ and a vector field $F$ on $U_z$ is given by the following formulas:
	\begin{enumerate}[label={\rm(\alph*)}]
		\item if $ C = \{z\}$, i.e.,  $z$ is an isolated critical point of $f$, then 
		\begin{equation}\label{eq:F-neaz-iso-zero}
			F_z = (f_z)'_y\frac{\partial}{\partial x} - (f_z)'_x\frac{\partial}{\partial y},
		\end{equation} 
		where $f_z$ is a local representation on $f$ on $U_z$ by a homogeneous polynomial without multiple factors  from {\rm (2a) of Definition \ref{def:class-F}};
		\item if $C$ is a critical circle of $f$, so $z\in C$ is a non-isolated critical point of $f$, then 
		\begin{equation}\label{eq:F_z-near-C-goal}
			F_z = \begin{cases}
				\frac{\partial}{\partial x}, & \text{if $C$ is a non-extremal critical circle of $f$},\\
				\mu_z(x,y) y \frac{\partial}{\partial x},  &\text{if $C$ is an extremal critical circle of $f$},
			\end{cases}
		\end{equation}
		for some  positive smooth  function $\mu_z:U_z\to \bbR$.
	\end{enumerate}
\end{enumerate}
\end{proposition}
\begin{definition}
A vector field $F$ as in {\rm Proposition~\ref{prop:vector-field-f}} will be called a Hamiltonian-like vector field\footnote{This definition provides a natural extension of the notion of a Hamiltonian-like vector field to functions with isolated singularities, which appears in \cite{Maksymenko:2021:review, Maksymenko:DefFuncI:2014}.}  of $f$, or for simplicity, an {\it$H$-field} of  $f$. 
\end{definition}
Proposition \ref{prop:vector-field-f} will be proved in Section \ref{sec:existence}, but now
we make several remarks.
Observe that a vector field $F$ has ``fewer'' singularities than the Hamiltonian field $X_f$. In particular, it has no zeros corresponding to the critical points of $f$ lying on the non-extremal critical circles of $f$. Explicit and simple local expressions for $F$ can be provided in neighborhoods of its zeros. Furthermore, in neighborhoods of the non-isolated critical points of $f$, the vector field $F$ is ``less  \textr{degenerate}'' than the Hamiltonian field $X_f$.
Throughout the paper, we will work with $F$ instead of $X_f$; this choice is motivated by Lemma~\ref{lm:hom-type-Did} and Proposition~\ref{prop:shift-functions-exist}.

\subsection{Period functions of $H$-like fields near critical circles}
Let $M$ be a smooth, connected, compact, and  \textr{orientable} surface, let $f$ be a function from $\F$, $F$ be an $H$-like field of $f$ with the flow $\bfF$. Let also $C$ be an extremal circle of $f$ and $Q$ be its saturated neighborhood.  So $Q$ is a cylinder such that $Q\setminus C$ contains no critical points of $f$.  Then each point $z\in Q\setminus C$ is periodic with respect to the flow $\bfF$ on $Q\setminus C$. Denote by $\theta:Q\setminus C\to \bbR$ a period function of $\bfF$ on $Q\setminus C$. It is known that $\theta$ is smooth on its domain.
The following result is standard, and we include its proof for completeness.
\begin{lemma}\label{lm:period}
\textr{The} period function $\theta:Q\setminus C\to \bbR$ of $\bfF$ on $Q\setminus C$ cannot be extended continuously  on $C$.
\end{lemma} 
\begin{proof}
This statement can be checked locally. Indeed, 	
let $z$ be a point at $C$. Then by definition of an $H$-field, there exists a local chart $(U,(x,y))$ near $z$ such that $z = (0,0)$, $\{p\in U\,|\,y(p) = 0 \}\subset C$ and
 $F = y\mu(x,y)\frac{\partial}{\partial x}$ on $U,$ where $\mu:U\to \bbR$ is a positive smooth function.  Recall that each point from $C\cap U$ is a zero of $F$; see \textsection\ref{sec:hor-lines}.

Let us fix $a,b> 0$  such that  $V = [0,a]\times [0,b]\subset U$.  We put $\max_{z\in V}\mu(z) = K>0$.
Since each point $(x,y)\in V$ with $y\neq 0$  is not a zero of $F$, it follows that  for any $y\in (0,b]$ there exists a unique number $t_y>0$ such that $\bfF_{t_y}(0,y) = (a,y)$. 
Integrating $F$ along the trajectory passing through $(0,y)$, $y\in (0,b]$, we get:
$$
a = x(t_y) = y\int_{0}^{t_y}\mu(x(\tau), y)d\tau\leq y\int_{0}^{t_y}Kd\tau = Kyt_y
$$
from which we immediately obtain
\begin{equation}\label{eq:t_y}
t_y\geq \frac{a}{Ky}\longrightarrow \infty\quad\text{as } y\to 0.
\end{equation}
Thus, the ``time'' required to move the point $(0,y)$ to $(a,y)$ along the trajectory of $\bfF$ tends to $\infty$ as $y\to 0.$ Consequently, \textr{the} period function $\theta:Q\setminus C\to \bbR$ of $\bfF$ cannot be extended to a continuous function $Q\to \bbR$.
\end{proof}

\section{Existence of $H$-fields for functions from $\F$}\label{sec:existence}
The aim of this section is to give the proof of Proposition \ref{prop:vector-field-f}.
Let $M$ be a smooth,  compact, connected, and \textr{orientable} surface, and let  $f$ be a function from $\F$ with the set of critical points $\Crit$.

In the next lemma, we prove the existence of the vector field  from  Proposition \ref{prop:vector-field-f} only in some neighborhoods of the critical circles of $f$.
\begin{lemma}\label{lm:vector-field-FC}
	Let $M$ be a smooth, compact, connected and \textr{orientable} surface, let $f$ be a function from $\F$ and $C$ be a critical circle of $f$. Then there exists a saturated neighborhood $Q_C$ of $C$ and a vector field $F_C$ on $\overline{Q_C}$ such that 
	\begin{enumerate}
		\item $f$ is constant along $F_C$ on $\overline{Q_C}$,
		\item for any $z\in C$ there exists a local chart $(U,(x,y))$ near $z$ such that $F_C$ on $U$ has the form
		\begin{align}
			F_C &= \frac{\partial}{\partial x},&& \text{if $C$ is non-extremal},\label{eq:FC-non-extremal}\\
			F_C &= \mu_z(x,y)y\frac{\partial }{\partial x}, && \text{if $C$ is extremal}\label{eq:FC-extremal},
		\end{align}
		for some smooth and positive function $\mu_z:U\to \bbR$.
	\end{enumerate}
\end{lemma}
\begin{proof}
Let also $z\in C$ be a non-isolated critical point of $f$. By (2.b) of Definition \ref{def:class-F} there exist a local coordinate system $(U_z, (\tilde{x},\tilde{y}))$ near $z$ such that a local representation of $f$ on $U_z$ is $f_z(\tilde{x},\tilde{y}) = a \tilde{y}^{n_C}$ for  some $n_C\geq 2$ depending on $C$, where $a \in \{ \pm 1\}$.

Let $\omega$ be a symplectic form on $M$.
 A Hamiltonian vector field $X_{f_z}$ of $f_z$ on $U_z$ has the form \eqref{eq:Ham-near-C}. By changing coordinate system $(x,y)= (a\tilde{x},\tilde{y})$ on $U_z$, we obtain that $X_f$ on $U_z$ has the form:
\begin{equation}\label{eq:Ham-chang-coord}
	X_f = \frac{n_Cy^{n_C-1}}{{g}(a x,y)}\frac{\partial}{\partial x},
\end{equation}
where $g:U_z\to \bbR$ is some positive and smooth function, see \textsection\ref{ssec:ham-function-F}.

The collection of sets $\{U_z\,|\,z\in C \}$ forms an open cover of $C$. From compactness of $C$, we can choose a finite ``chain-like'' subcover of $\{U_z\}_{z\in C}$. Namely,  
there exist $N\in \bbN$ and a subset $\{z_0, z_1,\ldots, z_{N-1}\}\subset C$ such that
\begin{enumerate}
	\item $\{U_{z_i}\}_{i = 0}^{N-1}$ is an open cover of $C$,
	\item $U_{z_i}\cap U_{z_j}\neq \varnothing$ iff $j = i\pm 1\;\mathrm{mod}\, N$.  
\end{enumerate}
{Choose a coordinate system $\phi_i = (x_i, y_i):U_{z_i}\to\mathbb{R}^2$ on $U_{z_i}$ with $\phi_i(C\cap U_{z_i})\subset \{y = 0\}$} such that  $X_f$ on $U_{z_i}$ has the form \eqref{eq:Ham-chang-coord} and
define a vector field $F_{z_i}$ on $U_{z_i}$ by the following formulas
\begin{align}
	F_{z_i} &= \frac{\partial}{\partial x_i}, &(\text{if $C$ is non-extremal, i.e., $n_C$ is odd}), \label{eq:F_z-C-non-extr}\\
	F_{z_i} &= y_i\frac{\partial}{\partial x_i}, &(\text{if $C$ is extremal, i.e., $n_C$ is even}). \label{eq:F_z-C-extr}
\end{align}
Note that the flows of vector fields $F_{z_i}$ and $X_f$ define the same foliation on $U_{z_i}\setminus C$.
If $C$ is non-extremal critical circle of $f$, then $U_z\cap C$ is a regular trajectory for $F_{z_i}$. 
From Eq.~\eqref{eq:Ham-chang-coord} it is easy to see that vector fields $F_{z_i}$  and $X_f$ are codirectional on $U_{z_i}\setminus C$.

Let $Q_C$ be a saturated neighborhood of $C$ such that $\overline{Q_C}\subset \bigcup_{i = 0}^{N-1} U_{z_i}.$ We will use the partition of unity to define a  vector field $F_C$ on $\overline{Q}_C$ mentioned above.
We put $Q_i = Q_C\cap U_{z_i}$ and $V_i = Q_i\setminus (U_{z_{i-1}}\cup U_{z_{i+1}})$, where all indexes are taken modulo $N$. It is easy to see that there exists a family smooth functions $\{ \rho_i:\overline{Q_C}\to \bbR\,|\, \rho_i\geq 0,\ i =0,1,\ldots, N-1 \}$
which satisfies the following conditions:
\begin{itemize}
	\item $\rho_i = 1$ on $\overline{V_i}$,
	\item $\rho_i = 0$ on $\overline{Q}\setminus \overline{Q_{i}}$,
	\item $\rho_i(z) + \rho_{i+1}(z) = 1$ for $z\in \overline{Q_i}\cap \overline{Q_{i+1}}$,
\end{itemize}
where  indexes are taken modulo $N$. Define a vector field $F_C$ on $\overline{Q_C}$ by the formula
\begin{equation}\label{eq:F_C-formula}
	F_C =\begin{cases}
		F_{z_i}, &\text{on }\overline{V_i},\\
		\rho_iF_{z_i}+\rho_{i+1} F_{z_{i+1}}, &\text{on } \overline{Q_i\cap Q_{i+1}}.
	\end{cases}
\end{equation}

Clearly, $F_C$ is codirectional with $X_f$ on $\overline{Q_C}\setminus C$ and $F_C(f|_{\overline{Q_C}}) = 0$. Therefore $f$ is constant along $F_C$ on $\overline{Q_C}$.

If $C$ is a non-extremal critical circle of $f$, then $F_C$ has no zeros in $\overline{Q_C}$. Thus,  for each $z\in C$ there exists a coordinate chart $(U, (x,y))$ near $z$ such that $F_C$ on $U$ has the form \eqref{eq:FC-non-extremal}.

{If $C$ is an extremal critical circle of $f$, then on $\overline{V_i}$ a vector field $F_C$ has the form $F_C = F_{z_i} = y_i\frac{\partial}{\partial x_i}$. It remains to prove that a vector field 
$\rho_iF_{z_i}+\rho_{i+1} F_{z_{i+1}}$ as above has the form $\eta_{i} F_{z_{i}}$ for some positive smooth function $\eta_{i}:\overline{Q_i\cap Q_{i+1}}\to \bbR$. This is a consequence of  Lemma \ref{lm:technical-lema}.

Indeed,  vector fields $F_{z_i}$ and $F_{z_{i+1}}$ have the same trajectories on $U_{z_i}\cap U_{z_{i+1}}$, and they are codirectional on $(U_{z_i}\cap U_{z_{i+1}})\setminus C$. Then there exist a ``special '' coordinate change, i.e., an orientation-preserving diffeomorphism 
$h:\phi_i(\overline{Q_i\cap Q_{i+1}})\to \phi_{i+1}(\overline{Q_i\cap Q_{i+1}})$ which satisfies conditions of the Lemma \ref{lm:technical-lema}, i.e., $h$ maps  trajectories of $(\phi_i)_*F_{z_i}$ to  trajectories of $(\phi_{i+1})_*F_{z_{i+1}}$ and preserve orientations of their regular trajectories. Thus, by (2) of Lemma \ref{lm:technical-lema}, on $\overline{Q_i\cap Q_{i+1}}$ a vector field
$\rho_iF_{z_i}+\rho_{i+1} F_{z_{i+1}}$ has the form $\eta_{i} F_{z_{i}}$ for some positive smooth function $\eta_{i}:\overline{Q_i\cap Q_{i+1}}\to \bbR$. This completes the proof of  Lemma \ref{lm:vector-field-FC}.}
\end{proof}

\subsection{Proof of Proposition \ref{prop:vector-field-f}}
Let $M$ be a smooth,  compact, connected, and \textr{orientable} surface, $\omega:TM\times TM\to\bbR$ be  a symplectic form on $M$, and let $X_f$  be a Hamiltonian vector field of $f$ with respect to $\omega,$ see \textsection\ref{ssec:Ham-skew}. A vector field $X_f$ satisfies  only (1) and (2) of Proposition \ref{prop:vector-field-f}. 

To obtain a vector field $F$, 
we define vector fields $F_w$  and $F_C$  in some neighborhoods of isolated critical points and critical circles of $f$ respectively by ``modifying'' $X_f$, and then using the partition of unity, ``replace'' $X_f$ on these neighborhoods by $F_w$ and $F_C$.

(a) {\it Vector field $F_w$ near isolated critical points of $f$.} Let $w$ be an isolated zero of $X_f$.  It follows from (2.a) of Definition \ref{def:class-F} that there exists a local coordinate system $(U_w, (x,y))$ near $w$ such that $f$ on $U_w$ is a homogeneous polynomial without multiple factors $f_w:\bbR^2\to \bbR$ of $\deg f_w> 1$ and so 
$X_f$ on $U_z$ has the form \eqref{eq:ham-vf-near-iso}.

Define a vector field $F_w$ on $U_w$ by the formula \eqref{eq:F-neaz-iso-zero}. 
A vector field $F_w$ and has the same trajectories as $X_f$ on $U_w$, $F_w$ and $X_f$ are codirectional on $U_w\setminus \{w\}$, and $F_w(f) = 0$ on $U_w,$ i.e., $F_w$ is tangent to leaves of $f$ on $U_w$. 

(b) {\it Vector field $F_C$ near  critical circles of $f$}. 
For each $C\in \CritC$ there exist a saturated neighborhood $Q_C$ of $C$ and a vector field $F_C$ on $\overline{Q_C}$ as in Lemma \ref{lm:vector-field-FC}.

(c) {\it Replacing $X_f$.} 
Using the partition of unity, we  modify the vector field $X_f$ by ``replacing'' it near the isolated critical points of $f$ with the vector fields $F_w$ as in (a), and near critical circles of $f$ with $F_C$ as in (b).
The resulting vector field $F$ on $M$ obviously satisfies conditions (1)--(3) of Proposition \ref{prop:vector-field-f}.\qed

 \section{Homotopy type of $\mathcal{D}_{\id}(\bfF)$ for  $H$-like fields}\label{sec:Hom-type-DidF}
 Let $M$ be a smooth, connected, compact, and  \textr{orientable} surface, let
 $f$ be a function from $\F$, and $F$ be a $H$-field of $f$. We  say that {\it $f$ has  property} (T) if at least one of the following conditions holds:
\begin{itemize}
	\item $f$ has a saddle (in this case $F$ has at least one non-closed trajectory),
	\item $f$ has  an isolated \textr{degenerate} local extreme (so $F$ has an isolated zero at which $1$-jet of $F$, i.e., a linear part of $F$, vanishes).
\end{itemize}

S.~Maksymenko \cite{Maksymenko:OsakaJM:2011} studied the group $\mathcal{D}_{\id}(\mathbf{Y})$ for a large class of flows $\{\mathbf{Y}: M \times \mathbb{R} \to M\}$ generated by vector fields on manifolds.
He proved that, under certain topological conditions at regular points and specific analytic conditions at singular points of $\mathbf{Y}$, every diffeomorphism from $\mathcal{D}_{\id}(\mathbf{Y})$ admits a shift function with respect to the flow $\mathbf{Y}$. This, in turn, made it possible to determine the homotopy type of $\mathcal{D}_{\id}(\mathbf{Y})$ for such flows.

The analytic conditions mentioned above specify the corresponding ``local normal forms'' of the vector field in neighborhoods of its zeros.
It should be noted that both the Hamiltonian vector field $X_f$ and an $H$-field $F$ associated with a function $f \in \F$ satisfy these conditions at regular points. However, unlike $F$, the Hamiltonian field $X_f$ generally fails to satisfy needed conditions at its singular points (see   Lemma \ref{lm:example}  for a local problem of the extending of  shift functions to singularities of the flows).

This provides the primary motivation for working with $F$ instead of $X_f$ in our setting. Theorem~3.5 from \cite{Maksymenko:OsakaJM:2011}, when applied to the $H$-fields of functions from $\F$, can be stated as follows:
\begin{lemma}[Theorem 3.5 \cite{Maksymenko:OsakaJM:2011}]\label{lm:hom-type-Did}
	Let $M$ be a smooth, connected, compact, and  \textr{orientable} surface, 
	let $f$ be a function from $\F$, $F$ be its $H$-like field with the flow $\mathbf{F}.$  
	Denote by $\Gamma_F^+$ the following convex subset of $C^{\infty}(M,\bbR)$
	$$
	\Gamma_F^+ = \{ \alpha\in C^{\infty}(M,\bbR)\,|\, 1+F\alpha>0 \}.
	$$
	Then the following holds true:
	\begin{enumerate}
		\item $\phi_{\bfF}(\Gamma^+_F) = \mathcal{D}_{\id}(\bfF),$ so each $h\in\mathcal{D}_{\id}(\bfF)$ has a shift function  $\alpha$ on $M$,
		\item  the restriction of a shift map 
		$$
		\phi_{\bfF}|_{\Gamma^+_F}:\Gamma^+_F\to \mathcal{D}_{\id}(\bfF)
		$$
		is either a homeomorphism  or a $\bbZ$-covering map. So $\mathcal{D}_{\id}(\bfF)$ is either contractible or has the homotopy type of $S^1$.
	\end{enumerate}
	In particular,
	if $f$ has either a property {\rm(T)} or $f$ has an extremal circle, then  $h\in\mathcal{D}_{\id}(\bfF)$ has a unique shift function $\alpha$ on $M$, \textr{$\phi_{\bfF}|_{\Gamma^+_F}$ is a homeomorphism, so 
	$\mathcal{D}_{\id}(\bfF)$ is contractible. }
	Otherwise, each point  $x\in M$ is periodic with respect to $\bfF$ and for each $n\in \bbZ$, the function $\alpha+n\theta$ is a shift function for $h$, where $\theta:M\to \bbR$ is a period function of $\bfF$.
\end{lemma}
\begin{proof}
	The only statement that is not completely covered in Theorem 3.5 \cite{Maksymenko:OsakaJM:2011} is the fact that if $f$ has at least one extremal circle, then $\mathcal{D}_{\id}(\bfF)$ is contractible. It follows from  Corollary 3.3. \cite{Maks:reparam-sh-map}.
	Here we present a simple proof of this fact.

	Assume that the converse is true, i.e, $f$ has at least one extremal circle and $\mathcal{D}_{\id}(\bfF)$ is homotopy equivalent to $S^1$. Then
	 by Corollary 3.3. \cite{Maks:reparam-sh-map}, this is equivalent to the fact that the vector field $G = \theta F$, where $\theta:M\to (0,\infty)$ is a period function of $F$ on $M$, yields a smooth circle action $\mathbf{G}:M\times\bbR\to M$, $\mathbf{G}(x,1) = x$ for each $x\in M$, where $\mathbf{G}$  is the flow of $G$.
	 This is not possible, since  for any regular neighborhood $Q$ of $C$ there exists a sequence of points $\{x_{\ell}\}_{\ell\geq 1}\subset Q\setminus C$ converging to $z\in C$ with $\lim\limits_{\ell\to \infty}\theta(x_{\ell}) = \infty$, see Lemma \ref{lm:period}.
	Therefore, our assumption is not true. 
\end{proof}

\section{Shift functions for diffeomorphisms from $\mathcal{S}_{\id}(f)$}\label{sec:shift-diff}
Let $M$ be a smooth, connected, compact, and  \textr{orientable} surface, let $f$ be a function from $\F$ with the set of extremal circles $E_f$ (possible empty), $F$ be an $H$-field of $f$ with the flow $\bfF$.
\subsection{Groups $\mathcal{S}_{\id}(f,E_f)$ and $\mathcal{G}(f,E_f)$}\label{ssec:groupsSidE}
{Denote by $\mathcal{D}(M,E_f)$ the group of diffeomorphisms of $M$ fixed on $E_f$. 
We put 
\begin{equation}\label{eq:SfEf}
\mathcal{S}(f,E_f) = \mathcal{S}(f)\cap \mathcal{D}(M,E_f)
\end{equation}
and denote by $\mathcal{S}_{\id}(f,E_f)$ a \textr{path} component of $\mathcal{S}(f,E_f)$ containing $\id_M$. So, if $h\in \mathcal{S}_{\id}(f,E_f)$, then there exists an isotopy $H:M\times [0,1]\to M$ such that for each $t\in [0,1]$ a diffeomorphism $H_t = H(-,t):M\to M$ satisfies the following properties: $f\circ H_t = f$, $H_0 = \id_M$, $H_1 = h$ and $H_t|_{E_f} = \id_{E_f}.$

We also set
\begin{equation}\label{eq:G}
\mathcal{G}(f, E_f) =  \Sid\cap \mathcal{D}(M,E_f).
\end{equation}
 So, if $h\in \mathcal{G}(f,E_f)$, then there exists an isotopy $H':M\times [0,1]\to M$ such that for each $t\in [0,1]$ a diffeomorphism $H'_t = H'(-,t):M\to M$ satisfies that following conditions: $f\circ H'_t = f$, $H'_0 = \id_M$, $H'_1 = h$ and $H'_1|_{E_f} = h|_{E_f} = \id_{E_f}$. So, in this case   $H'_t$ can ``move'' extremal circles on itself, that is $H'_t(C) \subset C$ for each $C\in E_f$, $t\in (0,1)$.}

It is easy to see that, in general, the group $\mathcal{G}(f,E_f)$ is not connected (the group $\pi_0\mathcal{G}(f,E_f)$ of \textr{path} components of $\mathcal{G}(f,E_f)$ will be studied in Section \ref{sec:pi0G}). Denote by  $\mathcal{G}_{\id}(f,E_f)$  a \textr{path} component $\mathcal{G}(f,E_f)$ which contains of $\id_M$. 

{
\begin{lemma}\label{lm:inclusions-stab}
	For $f\in \F$, the groups  $\mathcal{G}_{\id}(f,E_f) $ and $\mathcal{S}_{\id}(f,E_f)$ coincide. Thus, the following inclusions hold
	\begin{equation}\label{eq:Sid-inclusions}
		\Sid\supset \mathcal{G}(f,E_f)\supset \mathcal{S}_{\id}(f,E_f).
	\end{equation}
	If $E_f = \varnothing$, then all  groups in \eqref{eq:Sid-inclusions} coincide.
\end{lemma}
\begin{proof}
Obviously that  $\mathcal{G}(f,E_f)$ is a subgroup of $\mathcal{S}(f,E_f)$; this immediately gives an inclusion $\mathcal{G}_{\id}(f,E_f)\subset \mathcal{S}_{\id}(f,E_f)$.  The last inclusion can be verified directly. Indeed, 
let $h\in\mathcal{G}_{\id}(f,E_f)$ be a diffeomorphism. Then there exists an isotopy $H':M\times[0,1]\to M$ as above which is fixed on $E_f$, i.e., $H'_t|_{E_f} = \id_{E_f}$ for all $t\in[0,1]$. 
This means that  $h$ belongs to $\mathcal{S}_{\id}(f,E_f)$. The inclusion $\mathcal{S}_{\id}(f,E_f)\subset \mathcal{G}_{\id}(f,E_f)$ is trivial. Thus $\mathcal{S}_{\id}(f,E_f)$ and $\mathcal{G}_{\id}(f,E_f)$ coincide.

	The rest of lemma is obvious and follows from Eq.~\eqref{eq:SfEf} and Eq.~\eqref{eq:G}.
\end{proof}
}

\subsection{Shift functions for diffeomorphisms from $\Sid$}
The following results are devoted to studying the cases when diffeomorphisms from the groups  in \eqref{eq:Sid-inclusions}  admit shift functions with respect to $\bfF$.
Here we only present  results needed for our proofs in Section \ref{sec:pi0G}.

The following lemma is known for functions from $\mathcal{F}(M,P)$ with the  only isolated singularities.
\begin{proposition}[cf.  Theorem 1.3\cite{Maksymenko:AGAG:2006}]\label{prop:shift-functions-exist-global}
	Let $M$ be a smooth, compact, connected, and \textr{orientable} surface, let $f$ be a function from $\F$ with the set of extremal circles $E_f$ {\rm(}possible empty{\rm)}, let $F$ be an $H$-field of $f$ with the flow $\bfF$. Then the following equality holds:
	\begin{equation}\label{eq:SidDid}
		\mathcal{S}_{\id}(f,E_f) = \mathcal{D}_{\id}(\bfF).
	\end{equation}
\end{proposition}
\begin{proof}
The equality \eqref{eq:SidDid} can be proved directly.
The inclusion $\mathcal{D}_{\id}(\bfF)\subset \mathcal{S}_{\id}(f,E_f)$  is trivial. Let us show that the converse inclusion is also true. 
Let $h\in \mathcal{S}_{\id}(f,E_f)$ be a diffeomorphism. Then there exists an isotopy
$H:M\times [0,1]\to M$ such that a diffeomorphism $H_t = H(-,t):M\to M$ satisfies the following conditions:
\begin{itemize}
	\item $H_t:M\to M$ preserves $f$, i.e, 
	\begin{equation}\label{eq:f-stab}
		f\circ H_t = f,
	\end{equation}
	\item  $H_0 = \id_M$, $H_1 = h$,
	\item $H_t|_{E_f} = \id_{E_f}$,
\end{itemize}
for all $t\in [0,1]$.
Therefore, we have $H_t(f^{-1}(c))\subset f^{-1}(c)$ for all $c\in P$. 

Applying the chain rule to Eq. \eqref{eq:f-stab}, we obtain:
\begin{equation}\label{eq:chain-rule}
	\textr{	(Tf)_z = (T(f\circ H_t))_z = (Tf)_{H_t(z)}\circ (TH_t)_{z}}
\end{equation}	
for each $z\in M$. Since $H_t$ is a diffeomorphism for each $t\in [0,1]$, it follows that \textr{$(TH_t)_z\neq 0$} for each $z\in M$. Therefore, if a point $z\in M$ is regular (critical) for $f$, then $H_t(z)$ is regular (critical) point of $f$ for all $t\in [0,1]$. 

As the consequence,
we have that $H_t(\gamma)\subset \gamma$, $t\in [0,1]$ if $\gamma$ is a closed
regular trajectory  or an isolated fixed point of $\bfF$.  In particular, if $\gamma$ is an isolated zero of $F$, then $H_t|_{\gamma} = \id|_{\gamma}$ for $t\in [0,1]$, which yields that $H_t(\gamma')\subset \gamma'$ for all $t\in[0,1]$, where $\gamma'$ is a non-closed trajectory of $F$, i.e., a non-closed trajectory which corresponds to a connected component of the complement to some saddle of $f$.
If $\gamma$ is a non-isolated zero of $F$, then $\gamma\in C$ for some extremal circle $C$ of $f$. By assumption, we have $H_t|_{E_f} = \id|_{E_f},$ so each non-isolated zero of $F$ is fixed under the action of $H_t$, $t\in [0,1].$  Thus, for each trajectory $\gamma$ of $F$ we have $H_t(\gamma)\subset \gamma$, therefore $h = H_1\in \mathcal{D}_{\id}(\bfF).$
\end{proof}
\begin{corollary}\label{cor:hom-type-SidE}
	Let $f$ be a function from $\F$.  If
	\begin{enumerate}[label={\rm (\alph*)}]
		\item $E_f = \varnothing,$ then  $\mathcal{S}_{\id}(f,E_f) = \mathcal{S}_{\id}(f,\varnothing) = \Sid$ is contractible if $f$ has property {\rm (T)}; it has a homotopy type of $S^1$ otherwise,
		\item $E_f\neq \varnothing$, then $\mathcal{S}_{\id}(f,E_f)$ is always contractible.
	\end{enumerate}
\end{corollary}

{
\begin{proposition}[cf.  Theorem 3.5	\cite{Maksymenko:OsakaJM:2011}, Proposition 2.6 \cite{Maksymenko:NKolyv:2009},  Lemma 6.5 \cite{Maksymenko:DefFuncI:2014}, Lemma 3.5 \cite{Maksymenko:AGAG:2006}]\label{prop:shift-functions-exist} 
Let $M$ be a smooth, compact, connected, and \textr{orientable} surface, let $f$ be a function from $\F$ with the set of extremal circles $E_f\neq \varnothing$, let $F$ be an $H$-field of $f$ with the flow $\bfF$. For an extremal critical circle $C\in E_f$ of $f$, denote by $Q$ a saturated neighborhood of $C$. Then 
for any pair $(X,A)$ from the list 
	\begin{enumerate}
		\item $(Q\setminus C, \mathcal{S}_{\id}(f))$,
		\item $(Q, \mathcal{G}(f,E_f))$,
		\item  $( M\setminus E_f,\Sid)$,
		\item $(M,\mathcal{S}_{\id}(f,E_f))$,
	\end{enumerate}
	the following holds true: for any $h\in A$, the restriction $h_X:=h|_{X}$ admits a shift function $a_X:X\to \bbR$ w.r.t $\bfF_{X\times\bbR}, $ i.e., $h_X = \bfF_{\alpha_X}$. 
	
	In particular, if $f|_X:X\to P$ has either the property {\rm(T)} or has an extremal circle {\rm(}only cases {\rm (2)} and {\rm (4)}{\rm)}, then a shift function $\alpha_X$ for $h_X$ is unique. Otherwise, each point $x\in X$ is periodic of some period $\theta_X(x)$ w.r.t $\bfF_{X\times \bbR}$ and for any $n\in \mathbb{Z}$, a function $\alpha_X+n\theta_X$ is a shift function for $h_X$ on $X$, where $\theta_X:X\to\bbR$ is the period function of $\bfF|_{X\times\bbR}$. 
\end{proposition}}

{
\begin{remark}
	{\rm 	
	 Let us observe a few facts about cases (1) and (2) of Proposition \ref{prop:shift-functions-exist}. 
	 
	 Let $C$ be a critical circle of $f:M\to P$ and $Q$ be its saturated neighborhood. Note that the flow $\bfF|_{(Q\setminus C)\times \bbR}$ is nonsingular and each point $x\in Q\setminus C$ is periodic with some period $\theta(x)$ w.r.t the flow $\bfF|_{(Q\setminus C)\times}$, where $\theta:Q\setminus C\to \bbR$ is the period function of $\bfF|_{(Q\setminus C)\times \bbR}.$ From the other hand, the flow $\bfF|_{Q\times\bbR}$ has singularities --- each point $z\in C$ is a fixed point for $\bfF|_{Q\times \bbR}$.
	 By Lemma \ref{lm:period}, $\theta$ cannot be extended continuously to $Q$.

	 (a) Let $H$ be a diffeomorphism from $\mathcal{G}(f,E_f).$ By (2) of Proposition \ref{prop:shift-functions-exist}, a diffeomorphism $h = H|_Q$ has a unique shift function $\alpha:Q\to \bbR$ on $Q$ w.r.t the flow $\bfF|_{Q\times\bbR}$. From the other hand, by (1) of Proposition \ref{prop:shift-functions-exist}, for the restriction $h|_{Q\setminus C}$ there exists a countable many shift functions  on $Q\setminus C$, w.r.t $\bfF|_{(Q\setminus C)\times\bbR}$. Namely, there exists smooth function $\alpha':Q\setminus C\to \bbR$ such that for each $n\in \bbZ$, a function $\alpha'+n\theta$ is a shift function for $h|_{Q\setminus C}.$
	 The restriction $\alpha|_{Q\setminus C}$ of the shift function $\alpha$ of $h$  to $Q\setminus C$ is a shift function of $h|_{Q\setminus C}$. Therefore, there exists a number $n_0\in \bbZ$ such that $\alpha|_{Q\setminus C} = \alpha'+n_0\theta.$ 	
	 
	 The statement (2) of Proposition  \ref{prop:shift-functions-exist} is ``complementary'' to the one formulated above. By (1), the diffeomorphism $h|_{Q\setminus C}$ admits countably many shift functions 
	 $$
	 \{\ldots, \alpha'-2\theta, \alpha'-\theta, \alpha', \alpha'+\theta, \alpha'+2\theta, \ldots\}.
	 $$
	 on $Q\setminus C$ as above, while by (2),  {\it the only one} of them  can be extended to a smooth shift function on $Q$ for $h|_{Q}$ w.r.t $\bfF_{Q\times\bbR}.$
	 
	 (b) Similar phenomena take place for cases (3) and (4), when $f$ does not have property (T).
}
\end{remark}
}

{
\begin{proof}[{\bf Proof of Proposition} {\rm\ref{prop:shift-functions-exist}}]
	 It should be noted that the existence of shift functions for diffeomorphisms from $A$ on $X$ follows from the observation that each $h\in A$ preserves the orbits of $\bfF$ on $X$	and it is isotopic to $\id_M$ together with Proposition \ref{prop:shift-functions-exist-global} and  some general known results on the existence and extension of shift functions.

	 {\bf Case {\rm (4)}. The pair $(X,A) =  (M,\mathcal{S}_{\id}(f,E_f))$}. The proof in this case follows from Eq.~\eqref{eq:SidDid} and Lemma \ref{lm:hom-type-Did}.
	
	{\bf Cases {\rm (1)} and {\rm (3)}. The pair $(X,A)$ is either $(Q\setminus C,\mathcal{S}_{\id}(f))$ or $(X,A) = (M\setminus E_f, \mathcal{S}_{\id}(f))$}. The proofs in these cases proceed along the same lines.
	
	Let $h$ be a diffeomorphism from $A$. Then there exists an isotopy $H:X\times[0,1]\to X$ such that for all $t\in [0,1]$ a diffeomorphism $H_t = H(-,t):X\to X$ satisfies the following conditions: $f\circ H_t = f$, $H_0 = \id|_{X},$ $H_1 = h|_{X}$. By Theorem 5.25 \cite{Maksymenko:Indiana:2010}, for an isotopy $\{H_t\,|\, t\in [0,1]\}$ there exists a family of shift functions $\{\alpha_t:X\to \bbR\,|\, t\in [0,1]\}$ with $\alpha_0 = 0$ and such that $H_t = \bfF_{\alpha_t}$ on $X$. Therefore, $h = H_1 = F_{\alpha_1}$ on $X$.
	
	If $f|_{X}:X\to P$ has property (T), then from  Lemma 6.5. \cite{Maksymenko:DefFuncI:2014} and Lemma 5.2 \cite{Maksymenko:NKolyv:2009} follows that such shift function is unique.  Otherwise, each point $x\in X$ is periodic and thus, for each $a\in\bbZ$, a function $\alpha - a\theta$ is a shift function for $h$ on $X$, where $\theta:X\to \bbR$ is a period function of $\bfF$ on $X$. In particular, when $X = Q\setminus C$, a diffeomorphism $h|_{X}$ always has  countable many shift functions.
	
	{\bf Case {\rm (2)}. The pair $(X,A) =  (Q,\mathcal{G}(f,E_f)$}.
	 This statement can  be proved by hand using Lemma \ref{lm:example}.
	 First of all, we will define a shift function for $h\in\mathcal{G}(f,E_f)$ on some neighborhood of $C$ and then we extend it to $Q$.
	 
	 {\it Existence of the shift function near critical points on $C$}.
	 For each $z\in C$ there exists a chart $(U_z, (x,y))$ near $z$ such that $F$ on $U_z$ has the form $\mu_z(x,y) y \frac{\partial}{\partial x}$ for some positive smooth function $\mu_z:U_z\to \bbR$; see Proposition \ref{prop:vector-field-f}.  Sets $\{U_z\,|\, z\in C\}$ define an open cover of $C$ and $\bigcup_{z\in C}U_z\subset Q$.
	 Since $h$  preserves trajectories of the flow $\bfF|_{U_z\times \bbR}$  and
	 orientation of its regular trajectories on $U_z$, $z\in C$, it follows from 
	 Lemma \ref{lm:example} that for  each $z\in C$ there exists a smooth function $\alpha_z:U_{z}\to \bbR$ such that $h|_{U_{z}} = \bfF_{\alpha_z}$ on $U_{z}$.
	We claim that such shift function form $h$ on $U_z$ is unique. Indeed, let $\beta_z:U_z\to \bbR$ be another shit function for $h$ on $U_z$, so $h = \bfF_{\alpha_z} = \bfF_{\beta_z}$ on $U_z.$ 
	
	{ Note that $Q\setminus C$ has two connected components $Q^-$ and $Q^+$. 
	Each point $z\in Q^{\pm}$ is  periodic with some period $\theta^{\pm}(z)$ w.r.t $\bfF_{(Q^{\pm})\times\bbR}$, where  $\theta^{\pm}:Q^{\pm}\to \bbR$ is a smooth period function of $\bfF|_{(Q^{\pm})\times\bbR}$ Then $\ker\phi_{Q^{\pm},\bfF} = \{n\theta^{\pm}\,|\, n\in\bbZ\}$ and  by (3) of Lemma \ref{lm:ker-shift}, we get 
	$$
		\alpha_z-\beta_z = n^{\pm}_z\theta^{\pm}\qquad\text{on}\qquad U_z\cap Q^{\pm}
	$$		
	for some $n_z^{\pm}\in\bbZ.$
	A function $\alpha_z-\beta_z$ is smooth on $U_z$, but there is a sequence of points $\{x^{\pm}_{\ell}\}\subset U_z\cap Q^{\pm}$ converging to some point at $C\cap U_z$ such that $\lim_{\ell\to \infty}\theta^{\pm}(x^{\pm}_{\ell}) = \infty$, see Lemma \ref{lm:period}. Thus, $n^{\pm}_z$  must be $0$, and therefore $\alpha_z = \beta_z$ on $U_z$.
}

	 Next we need to show that such shift functions match at the overlaps $U_z\cap U_w\neq \varnothing$.

	{\it Existence of shift function near $C$}.
	 Let $z,w\in C$ be  points such that $U_z\cap U_w\neq \varnothing$, $\alpha_z:U_z\to \bbR$ and $\alpha_w:U_w\to \mathbb{R}$ be shift functions of $h$ on $U_z$ and $U_w$ respectively as above. Then on $U_{z}\cap U_{w}$ we have $h = \mathbf{F}_{\alpha_{z}} = \bfF_{\alpha_w}$.
	  {Again, by (3) of Lemma \ref{lm:shift-main} we have 
	$$
	\alpha_z -\alpha_w = n^{\pm}_{zw}\theta^{\pm}\qquad\text{on}\qquad (U_{z}\cap U_{z_w})\cap Q^{\pm}
	$$	 
	for some $n_{zw}\in \mathbb{Z}$, where
	$\theta^{\pm} :Q^{\pm}\to \mathbb{R}$ is a period function for
	$\mathbf{F}$ on $Q^{\pm}$.
 	Similarly to the above, it is easy to see that $n^{\pm}_{zw}$ must be $0$, so  we get $\alpha_{z} = \alpha_w$ on $U_{z}\cap U_{w}$. This shows that  functions $\alpha_z$ and $\alpha_w$ for $h$ on $U_z$ and $U_w$ agree on overlaps $U_z\cap U_w$. 
 	Therefore, there exist a unique smooth function $\alpha:W\to \bbR$, where $W =  \bigcup_{z\in C}U_{z}$  such that $\alpha|_{U_{z}} = \alpha_z$ and  $h = \bfF_{\alpha}$ on $W$.
}
	  
	 {\it Extension of shift function on $Q$}.
	 By (1) of Proposition \ref{prop:shift-functions-exist}, there exists a smooth function $\alpha':Q\setminus C\to \mathbb{R}$ such that for any $n\in\bbZ$, a function $\alpha'+n\theta$ is a shift function for $h|_{Q\setminus C}$. Then there exists $n_0\in \bbZ$ such that $\alpha|_{W\setminus C}$  has the form $ \alpha'+n_0\theta$ on $W\setminus C$. Thus, a shift function $\alpha$ for $h$ on $W$ can be extended to a smooth shift function $\alpha:Q\to \bbR$ for $h$ on $Q$.
\end{proof}

}

\section{The fibration for $\Sid$}\label{sec:fibration-rho}
Let $M$ be a smooth, connected, compact and  \textr{orientable} surface, and let $f$ be a function from $\F$ with the set of extremal circles $E_f$. {\it Throughout this section, we assume that $E_f\neq \varnothing$}. 
Note that $E_f$ is a compact $1$-manifold. Denote by $\mathcal{D}(E_f)$ the group of diffeomorphisms of $E_f$ and by $\mathcal{D}_{\id}(E_f)$ a \textr{path} component of $\mathcal{D}(E_f)$ containing $\id_{E_f}$. 
{The following lemma is well-known.
\begin{lemma}[cf.~Proposition 4.2, \cite{Ghys:Enseign:2001}]\label{lm:DidE}
	The group $\mathcal{D}_{\id}(E_f)$ is locally contractible\footnote{By locally contractibility  we mean {\it strong} local contractibility. Recall that a space is strongly locally contractible if each point has a local base of contractible neighborhoods.} and has the homotopy type of a torus $T^{|E_f|}$.
\end{lemma}
\begin{proof}
	The group  $\mathcal{D}_{\id}(S^1)$ is known to be a Fr\'echet manifold, diffeomorphic to $S^1\times L,$ where 
	$$
	L = \Bigg\{f\in C^{\infty}(\mathbb{R}, \bbR)\,\Big|\, f\text{ is strictly positive and $2\pi$-periodic with }  \frac{1}{2\pi}\int_0^{2\pi} f(t)\,dt = 1 \Bigg\}.
	$$
	The space $L$ is an open and convex subspace of the space $C^{\infty}(\bbR, \bbR)$, thus $L$ is contractible and  locally contractible topological space. The contractibility of $L$ implies the homotopy equivalence $\mathcal{D}_{\id}(S^1)\sim S^1$, and the local contractibility of $L$ implies the local contractibility of $\mathcal{D}_{\id}(S^1)$.
	
	Assume that $E_f = \{C_1,C_2,\ldots, C_n\}$, $n\geq 1$. Note that curves from  $E_f$ are pairwise disjoint. Then $\mathcal{D}_{\id}(E_f)$ is naturally isomorphic to $\prod_{i = 1}^n\mathcal{D}_{\id}(C_i)$, which yields that $\mathcal{D}_{\id}(E_f)$ is locally contractible and $\mathcal{D}_{\id}(E_f)\sim (S^1)^n$.
\end{proof}
}

 O.~Khokhliuk and S.~Maksymenko \cite{KhohliykMaksymenko:Indag:2020} studied diffeomorphisms of surfaces which preserve the given foliation with singularities on $M$ with some ``nice'' behavior near critical leaves.  We formulate their main result adapted to our case.
\begin{theorem}[Theorem 8.2., Theorem 3.3.  \cite{KhohliykMaksymenko:Indag:2020}]\label{thm:fibration-rho}
	Let $f$ be a function from $\mathcal{F}(M,P)$ with   $E_f = \{C_1, \ldots, C_n\}$, $n\geq 1$.
	Then the map 
	\begin{equation}\label{eq:fibration-rho}
			\rho:\mathcal{S}(f)\to \mathcal{D}(E_f),\qquad\rho(h) = h|_{E_f}
	\end{equation}
is a locally-trivial fibration with the fiber $ \mathcal{S}(f,E_f)$. In particular, the restriction 
\begin{equation}\label{eq:fibration-rho0}
	\rho_0:\mathcal{S}_{\id}(f)\to \mathcal{D}_{\id}(E_f) =\prod_{i = 1}^n\mathcal{D}_{\id}(C_i), \qquad \rho_0(h) = (h|_{C_1}, h|_{C_2}, \ldots, h|_{C_n})
\end{equation}
is also a locally-trivial fibration with the fiber $\mathcal{G}(f,E_f) = \mathcal{D}(M, E_f)\cap \Sid$; see  {\rm Eq.~\eqref{eq:G}}.
\end{theorem}
A long sequence of homotopy groups of the fibration \eqref{eq:fibration-rho0} will be our main tool for the study of homotopy properties of $\Sid$ for $f\in \F$ with $E_f\neq \varnothing$.
The following lemmas are  consequences of Theorem \ref{thm:fibration-rho} and Corollary \ref{cor:hom-type-SidE}.

\begin{lemma}\label{lm:CW}
	Let $f$ be a function from $\F$ with $E_f = \{C_1, C_2, \ldots, C_n\}$, $n\geq 1$. Then 
	$\Sid$ has the homotopy type of a CW complex.
\end{lemma}
\begin{proof}
{(1) Firstly we claim that $\mathcal{G}(f,E_f)$ is a locally contractible topological group.

Let $F$ be an $H$-field of $f$ and $\bfF:M\times \bbR\to M$ be its flow. By Lemma \ref{lm:inclusions-stab} and Eq.~\eqref{eq:SidDid}, we have
$
\mathcal{G}_{\id}(f,E_f) = \mathcal{S}_{\id}(f,E_f) = \mathcal{D}_{\id}(\bfF).
$
By assumption we have $E_f\neq \varnothing$, so  from  Lemma \ref{lm:hom-type-Did} follows that  $\mathcal{G}_{\id}(f,E_f)$ is homeomorphic to a convex subset $\Gamma_F^+ = \{ \alpha\in C^{\infty}(M,\bbR)\,|\, 1+F\alpha>0 \}$ of $C^{\infty}(M,\bbR)$ which is  locally contractible. Since each path  component of $\mathcal{G}(f,E_f)$ is homeomorphic to $\mathcal{G}_{\id}(f,E_f)$, it follows that
 the group $\mathcal{G}(f,E_f)$ is  locally contractible.

(2) The group $\mathcal{S}_{\id}(f)$ inherits many topological properties of the group $\mathcal{D}(M)$. In particular, the group $\mathcal{S}_{\id}(f)$ is Hausdorff and paracompact. 
We claim that $\mathcal{S}_{\id}(f)$ is a locally contractible topological group.
Then by results of R.~Palais \cite{Palais:MemAMS:1960} and J.~Milnor \cite{Milnor:Trans:1959}, $\mathcal{S}_{\id}(f)$ has the homotopy type of a CW complex.

Local contractibility $\Sid$  follows from the properties of the fibration $\rho_0$, see \eqref{eq:fibration-rho0}. 
To be more precise, the fibration $\rho_0$ is locally trivial, i.e., $\rho_0$ admits local sections.
Let $h$ be a diffeomorphism from $\Sid$.
By Lemma \ref{lm:DidE},
for $\rho_0(h)\in \mathcal{D}_{\id}(E_f)$, there  exists an open and  contractible neighborhood $U$ such that $\rho_0^{-1}(U)\cong U\times \mathcal{G}(f, E_f)$.  Since $\mathcal{G}(f,E_f)$ is locally contractible (see above) it follows that  $ U\times \mathcal{G}(f, E_f)$ is locally contractible. This completes the proof of locally contractibility of  $\Sid$.}
\end{proof}

\begin{lemma}\label{lm:weak-Sid}
	Let  $f$ be a function from $\F$ with the set of extremal circles $E_f = \{C_1,\ldots, C_n\}$, $n\geq 1$. Then the following holds true:
	\begin{enumerate}
		\item ${\pi_q\Sid} = 0$ for $q\geq 2$,
		\item  there is a short exact sequence of homotopy groups for {$\pi_1\Sid$}:
		\begin{equation}\label{eq:non-triv-subseq}
			\xymatrix@R=0.5cm{
				1\ar[r]& {\pi_1\Sid} \ar[r]^{\rho_1}&\pi_1\mathcal{D}_{\id}(E_f) \ar@{=}[d] \ar[r]^{\partial_1}& \pi_0\mathcal{G}(f, E_f) \ar[r] & 1\\
				&& {\mathbb{Z}^n}		&		
			}
		\end{equation}
	\end{enumerate}
\end{lemma}
\begin{proof}
	There is a long exact sequence of homotopy groups of the fibration $\rho_0$ with base point $\id_M$:
	\begin{align}\label{eq:hom-group-seq-fibrat}
		\begin{split}
			\ldots \longrightarrow &\pi_{q+1}\mathcal{D}_{\id}(E_f) \longrightarrow \pi_{q} \mathcal{G}_{\id}(f, E_f)\longrightarrow\pi_q\Sid\longrightarrow \pi_q \mathcal{D}_{\id}(E_f) \longrightarrow\ldots\\
			&\ldots \longrightarrow\pi_1\mathcal{D}_{\id}(E_f)\longrightarrow \pi_0\mathcal{G}(f, E_f)\longrightarrow \pi_0\Sid\longrightarrow\pi_0\mathcal{D}_{\id}(E_f)\longrightarrow 1.
		\end{split}
	\end{align}
	First, we start with some ``replacements'' in the sequence \eqref{eq:hom-group-seq-fibrat}.
	The groups $\mathcal{S}_{\id}(f)$ and $\mathcal{D}_{\id}(E_f)$ are path-connected, thus $\pi_0\mathcal{S}_{\id}(f) = \pi_0\mathcal{D}_{\id}(E_f) = 1.$
	The identity path component $\mathcal{G}_{\id}(f, E_f)$ of $\mathcal{G}(f, E_f)$ is $\mathcal{S}_{\id}(f,E_f)$, so in \eqref{eq:hom-group-seq-fibrat} we replace $\pi_q\mathcal{G}_{\id}(f,E_f)$ by $\pi_q\mathcal{S}_{\id}(f,E_f)$ for $q\geq 1$.
	
	{ By Lemma \ref{lm:DidE},  the group $\mathcal{D}_{\id}(E_f) = \prod_{i = 1}^n\mathcal{D}_{\id}(C_i)$ has the homotopy type of an $n$-torus $T^n$, so }
	$$
	\pi_q\mathcal{D}_{\id}(E_f) = \begin{cases}
		\bbZ^n, & \text{if }q = 1,\\
		1, &\text{otherwise},
	\end{cases}
	$$
	which yields an isomorphism $\pi_q\mathcal{G}_{\id}(f, E_f) = \pi_q\mathcal{S}_{\id}(f,E_f)\cong \pi_q\Sid$ for $q\geq 2.$ 
	
	By assumption,  $E_f\neq\varnothing$, then from Corollary \ref{cor:hom-type-SidE} the group $\mathcal{S}_{\id}(f,E_f)$ is contractible. Therefore, $\pi_q\Sid \cong \pi_q\mathcal{S}_{\id}(f,E_f) = 0$ for $q\geq 2$, and
	the sequence  \eqref{eq:non-triv-subseq} is a non-trivial part of the tail of a long exact sequence \eqref{eq:hom-group-seq-fibrat}.
\end{proof}
By Lemma \ref{lm:CW}, Lemma \ref{lm:weak-Sid} and Whitehead theorem, the homotopy type of $\mathcal{S}_{\id}(f)$ depends only on the group  $\pi_1\mathcal{S}_{\id}(f)$, which is a part of  a short exact sequence \eqref{eq:non-triv-subseq}. The group $\pi_0\mathcal{G}(f, E_f)$ will be studied in Section \ref{sec:pi0G}.

\section{Description of $\pi_0\mathcal{G}(f,E_f)$}\label{sec:pi0G}
Let $M$ be a smooth, compact, connected, and \textr{orientable} surface, and let $f$ be a function from $\F$ with the set of extremal circles $E_f = \{C_1,C_2,\ldots, C_n\}$, $n\geq 1$. Let also $F$ be an $H$-field of $f$ with the flow $\bfF:M\times\bbR\to M$. The aim of this section is to give the description of the group $\pi_0\mathcal{G}(f, E_f)$, where $\mathcal{G}(f,E_f) = \mathcal{D}(M, E_f)\cap \Sid$, see \textsection\ref{ssec:groupsSidE}, and
our main result   is the following proposition.
\begin{proposition}\label{prop:pi0-G}
	The group $\pi_0\mathcal{G}(f, E_f)$ is a free abelian group. If
	$f$ satisfies {\rm (T)},  then $\pi_0\mathcal{G}(f, E_f)$ is isomorphic to $\bbZ^{|E_f|}$; otherwise, $\pi_0\mathcal{G}(f, E_f)$ is isomorphic to $\bbZ^{|E_f|-1}$.
\end{proposition}

\subsection{Proof of Proposition \ref{prop:pi0-G}}
 Let $Q_i$ be a saturated neighborhood of $C_i\in E_f$, $i = 1,2,\ldots, n$.  Note that the flow $\bfF|_{(Q_i\setminus C_i)\times\bbR}$ on $Q_i\setminus C_i$ has no fixed points.

Let $h$ be a diffeomorphism from $\mathcal{G}(f, E_f).$ 
By case (2) of Proposition \ref{prop:shift-functions-exist} there exists a unique smooth function $\beta_i:Q_i\to \bbR$ such that $h|_{Q_i} = \bfF_{\beta_i}$ on $Q_i$. 
Here the proof splits into two cases (A) and (B).

{\bf Case (A).} {\it Assume that $f$ has property} (T), i.e., $f$ has  either an isolated \textr{degenerate} local extreme, or a saddle. Then,  by case (3) of  Proposition \ref{prop:shift-functions-exist},  there exists a unique smooth function $\alpha:M\setminus E_f\to \bbR$ such that $h = \bfF_{\alpha}$ on $M\setminus E_f.$ Thus,  we have $h(x) = \bfF(x,\alpha(x)) = \bfF(x,\beta(x))$ for  $x\in Q_i\setminus C_i$. Since every point in $Q_i\setminus C_i$ is periodic with respect to $\bfF|_{(Q_i\setminus C_i)\times\bbR}$ it follows that there exists $c_i(h)\in\bbZ$ such that 
\begin{equation}\label{eq:c_i(h)-non-per}
\alpha-\beta_i = c_i(h)\theta_i\qquad \text{on }Q_i\setminus C_i,
\end{equation}
where $\theta_i:Q_i\setminus C_i\to \bbR$ is a period function of $\bfF|_{(Q_i\setminus C_i)\times\bbR}$.

\begin{lemma}\label{lm:zeta-alpha-unique}
	 A map $\zeta:\mathcal{G}(f, E_f)\to \bbZ^n$ defined by
	\begin{equation}\label{eq:map-ZZZn}
		\zeta(h) = (c_1(h), c_2(h), \ldots, c_n(h)),
	\end{equation}
where numbers $c_i(h)$ are given by {\rm Eq.~\eqref{eq:c_i(h)-non-per}},
 is a homomorphism with $\ker\zeta \cong \mathcal{S}_{\id}(f, E_f)$. 
\end{lemma}
\begin{proof}
	(1) First, we show that $\zeta$ is a homomorphism.
	Let $g$ be a diffeomorphism from $\mathcal{G}(f, E_f)$. Then, by Proposition  \ref{prop:shift-functions-exist}, there exists unique smooth functions $\gamma:M\setminus E_f\to \bbR$ and $\delta_i:Q_i\to \bbR$, $i=1,2,\ldots, n$ such that $g = \bfF_{\gamma}$ on $M\setminus E_f$ and $g = \bfF_{\delta_i}$ on $Q_i\to \bbR$.
	
	Assume that $\zeta(g) = (c_1(g), c_2(g), \ldots, c_n(g))\in\bbZ^n$,  where $c_i(g)$ is given by
	\begin{equation}\label{eq:gamma-delta}
		\gamma-\delta_i = c_i(g)\theta_i, \qquad \text{on }Q_i\setminus C_i.
	\end{equation}
	By (3) of Lemma \ref{lm:shift-main}, we obtain that $h\circ g = \bfF_{\sigma}$ on $M\setminus E_f$ and $h\circ g = \bfF_{\varkappa_i}$ on $Q_i$, where 
	$$
	\sigma = \gamma+\alpha\circ \bfF_{\gamma}\qquad \varkappa_i =\delta_i+ \beta_i\circ \bfF_{\delta_i}.
	$$
	 From the one hand,  $\zeta(h\circ g) = (c_1(h\circ g), c_2(h\circ g),\ldots,  c_n(h\circ g))$, i.e.,
	 \begin{equation}\label{eq:zeta-gen}
	 \sigma-\varkappa_i = c_i(h\circ g)\theta_i,\qquad\text{on }Q_i\setminus C_i,
	 \end{equation}
	and from the other hand on $Q_i\setminus C_i$ we have
	\begin{align*}
	\sigma-\varkappa_i &=\gamma+\alpha\circ \bfF_{\gamma}-\delta_i-\beta_i\circ \bfF_{\delta_i} \\
	&=  (\alpha\circ \bfF_{\gamma}-\beta_i\circ \bfF_{\delta_i}) +(\gamma-\delta_i ) & \\
	&=  (\alpha\circ \bfF_{\gamma}-\beta_i\circ \bfF_{\gamma-c_i(g)\theta_i}) + (\gamma-\delta_i )
	& (\text{from Eq. \eqref{eq:gamma-delta} }\delta_i = \gamma-c_i(g)\theta_i)\\
	&= ( \alpha\circ \bfF_{\gamma} - \beta_i\circ \bfF_{\gamma}) +(\gamma-\delta_i ) & (\text{$\theta_i$ is a period function on $\bfF$ on $Q_i\setminus C_i$} )\\
	&=  (\alpha-\beta_i)\circ \bfF_{\gamma} +(\gamma-\delta_i ) \\
	&= c_i(h)\theta_i\circ \bfF_{\gamma} +c_i(g)\theta_i \\
	&=c_i(h)\theta_i +c_i(g)\theta_i & (\text{since $\theta_i$ is constant on each trajectory of $\bfF$}) \\
	&= (c_i(h)+c_i(g))\theta_i.
	\end{align*}
	Then $c_i(h\circ g) = c_i(h)+c_i(g)$ for each $i = 1,2,\ldots, n$ which yields that $\zeta(h\circ g) = \zeta(h)+\zeta(g)$. Therefore,  $\zeta$ is a homomorphism.
	
	(2) It remains to prove that $\ker\zeta = \mathcal{S}_{\id}(f,E_f)$.
	 Let $h$ be a diffeomorphism from $\mathcal{G}(f,E_f)$, and let $\alpha$ and $\beta_i$, $i = 1,2,\ldots, n$ are shift functions of $h$ on $M\setminus E_f$ and $Q_i$ as above. 

	Assume that $h\in \ker\zeta.$ Then $c_i(h) = 0$ for all $i = 1,2,\ldots, n,$ and thus, $\alpha = \beta_i$ on $Q_i\setminus C_i$. A function $\beta_i$ is smooth on $Q_i$, $i = 1,2,\ldots, n$, so a function $\alpha$ on $M\setminus E_f$ can be smoothly extended to the set of critical circles $E_f$ of $f$; the resulting function $\tilde{\alpha}:M\to \bbR$ is a shift function for $h$ on $M$, i.e.,  $h = \bfF_{\tilde{\alpha}}$ on $M$. Then, by case (4) of Proposition \ref{prop:shift-functions-exist}, a diffeomorphism $h$ belongs $\mathcal{S}_{\id}(f,E_f)$.
		
	Assume that $h$ belongs to $\mathcal{S}_{\id}(f,E_f)$. Since $E_f \neq\varnothing$, it follows from (4) of Proposition \ref{prop:shift-functions-exist} that there exists a unique smooth function $\alpha:M\to \bbR$ such that $h = \bfF_{\alpha}$ on $M$.  By Eq. \eqref{eq:c_i(h)-non-per}, the following holds:
	$\alpha-\beta_i = c_i(h)\theta_i$ on $Q_i\setminus C_i$. Note that functions $\alpha,\beta_i$ are smooth on $Q_i$, and $\theta_i$ is smooth on $Q_i\setminus C_i$. But there exists a sequence $\{x_{\ell_i}|, {\ell_i}\geq 1\}$ of points from $Q_i\setminus C_i$ converging to some point $z\in C_i$  such that  $\lim_{\ell_i\to \infty}\theta(x_{\ell_i}) = \infty$, see Lemma \ref{lm:period}. Then $c_i(h)$ must be $0$ for all $i = 1,2,\ldots, n$, and therefore, $h$ belongs to $\ker\zeta.$	 
\end{proof}
By Lemma \ref{lm:zeta-alpha-unique}, we have $\mathcal{G}(f,E_f)/\ker \zeta\cong \bbZ^n$.
Finally,  using $\mathcal{G}_{\id}(f,E_f) = \mathcal{S}_{\id}(f,E_f)$, we get
$$
 \pi_0\mathcal{G}(f,E_f) =  \mathcal{G}(f,E_f)/\mathcal{G}_{\id}(f, E_f)  = \mathcal{G}(f, E_f)/\mathcal{S}_{\id}(f,E_f) \cong \mathcal{G}(f, E_f)/\ker\zeta\cong \bbZ^n,
$$
which ends the proof in this case.

{\bf Case (B).} {\it Assume that $f$ does not have property} (T), so all isolated critical points of $f$ are non-\textr{degenerate} (Morse) local extremes.
Then
each point $z\in M\setminus E_f$ is periodic with respect to $\bfF$ on $M\setminus E_f$ and,
 by case (3) of Proposition \ref{prop:shift-functions-exist},
a shift function for $h$ on $M\setminus E_f$ is not unique. Namely, 
there exists a smooth function $\alpha:M\setminus E_f\to \bbR$ such that for any $a\in \bbZ$, the function $\alpha-a\theta$ is a shift function for $h$ on $M\setminus E_f$, i.e., $h = \bfF_{\alpha-a\theta}$, where $\theta:M\setminus E_f\to \bbR$ is a period function of $\bfF$ on $M\setminus E_f$.

For $a\in \bbZ$ and  a shift function $\alpha-a\theta$ for $h$ on $M\setminus E_f$,
the following holds: $h(x) = \bfF(x,\alpha(x) -a\theta(x)) = \bfF(x,\beta_i(x))$ for $x\in Q_i\setminus C_i$. Since each point in $Q_i\setminus C_i$ is periodic, it follows  that
there exists a unique $c_i(h)\in \bbZ$ such that $\alpha -a\theta-\beta_i = c_i(h)\theta_i$ on $Q_i\setminus C_i$, $i = 1,2,\ldots, n$, where $\theta_i:Q_i\setminus C_i\to \bbR$ is a period function of $\bfF|_{(Q_i\setminus C_i)\times\bbR}$.
 Note that  $\theta|_{Q_i\setminus C_i} = \theta_i$ which yields
\begin{equation}\label{eq:c_i(h)-per}
	\alpha -\beta_i = (c_i(h)+a)\theta_i\qquad \text{on }Q_i\setminus C_i.
\end{equation}

Compared to the case (A), numbers in RHS in Eq. \eqref{eq:c_i(h)-per} depend on the choice of a shift function  for $h$, and
thus, they are not well-defined for $h$. However, ``the vector  $(c_1(h),c_2(h),\ldots, c_n(h))$ is well-defined for $h$ up to  sums with  constant vectors'' $(a,a,\ldots, a)\in\bbZ^n$.
To make this precise, we recall the following construction.

{\it Short aside. Quotient by diagonal.}
Let 
$$
\Delta = \{\underbrace{(a,a,\ldots, a)}_{n}\,|\,a\in\bbZ\}
$$ be a subgroup of $\bbZ^n$.
Consider the following equivalence relation:
two vectors $c = (c_1,c_2,\ldots, c_n)\in\bbZ^n$ and $d = (d_1,d_2,\ldots, d_n)\in\bbZ^n$ are equivalent  $c\sim d$ iff $c-d\in \Delta,$ i.e., there exists $a\in \bbZ$ such that 
$$
(d_1,d_2,\ldots, d_n) = (c_1+a,c_2+a,\ldots, c_n+a).
$$
Denote by $[c] = c+\Delta = \{c+(a,a,\ldots, a)\,|\,a\in \bbZ\}$ the equivalence class of $c\in \bbZ^n$. Let $\bbZ^n/\Delta = \{ [c] \,|\, c\in \bbZ^n\}$ be a quotient group of $\bbZ^n$ by $\sim$ with the addition defined by $[c] + [d]= [c+d]$ and let 
$$
p:\bbZ^n\to \bbZ^n/\Delta,\qquad p(c) = [c] = c+\Delta,
$$
be a canonical projection. Therefore, $c\sim d$ iff $p(c) = p(d).$
 It is also known that $\bbZ^n/\Delta$ is isomorphic to $\bbZ^{n-1}$ via an isomorphism $[c_1,c_2,\ldots, c_n]\mapsto (c_2-c_1, c_3-c_1, \ldots, c_n-c_1)$.

 Returning to the proof, one can state the following result, similar  to Lemma \ref{lm:zeta-alpha-unique}.
\begin{lemma}\label{lm:zeta-alpha-not-unique}
	A map $\psi:\mathcal{G}(f, E_f)\to \bbZ^n/\Delta$ defined by
	\begin{equation}\label{eq:map-ZZZn-Delta}
		\psi(h) = [c_1(h), c_2(h), \ldots, c_n(h)]\in \bbZ^n/\Delta,
	\end{equation}
where numbers $c_i(h)$ are given by {\rm Eq.~\eqref{eq:c_i(h)-per}},
	is a homomorphism with $\ker\psi = \mathcal{S}_{\id}(f, E_f)$. 
\end{lemma}
\begin{proof}
This result can be proved similarly to Lemma \ref{lm:zeta-alpha-unique}.
We begin by verifying that $\psi$ is a homomorphism.	

(1) 
 Let $g$ be diffeomorphism from $\mathcal{G}(f,E_f)$ with shift functions $\gamma-b\theta:M\setminus E_f\to \bbR$, where $b\in \bbZ$ and $\delta_i: Q_i\to \bbR$ on $Q_i$, $i = 1,2,\ldots, n.$ Then by  Eq. \eqref{eq:c_i(h)-per}, there exists a unique $c_i(g)\in \bbZ$ such that 
\begin{equation}\label{eq:diff-on-Q-C}
\gamma-\delta_i = (c_i(g)+b)\theta_i, \qquad \text{on }Q_i\setminus C_i,
\end{equation}
and thus, $\psi(g) = [c_1(g), c_2(g), \ldots, c_n(g)]\in\bbZ^n/\Delta.$

Since $h$ and $g$ has countable number of shift functions on $M\setminus E_f$, it follows that $h\circ g$ also has countable number of shift functions on $M\setminus E_f$. One can show, using (3) of Lemma \ref{lm:shift-main} that for each $d\in \bbZ$, a function $\sigma -d\theta$, where $\sigma = \gamma+\alpha\circ \bfF_{\gamma}$, is a shift function for $h\circ g$ on $M\setminus E_f$, i.e. $h = \bfF_{\sigma-d\theta}$ on $M\setminus E_f$.  In particular, for fixed $a,b\in\bbZ$, we have $d = a+b$.

Since $h$ and $g$ have unique shift functions on $Q_i$, it follows that
a diffeomorphism $h\circ g$ has a unique shift function $\varkappa_i$ on $Q_i$ and, again by (3) of Lemma \ref{lm:shift-main}, $\varkappa_i = \delta_i+\beta_i\circ \bfF_{\delta_i}.$

From the one hand, by Eq. \eqref{eq:c_i(h)-per},  we get
$$
\sigma-\varkappa_i = (c_i(h\circ g) +d)\theta_i,\quad\text{on } Q_i\setminus C_i,
$$
so $\psi(h\circ g) = [c_1(h\circ g), c_2(h\circ g), \ldots, c_n(h\circ g)],$
 and from the other hand on $Q_i\setminus C_i$ we have:
\begin{align*}
	\sigma-\varkappa_i &=\gamma+\alpha\circ \bfF_{\gamma}-\delta_i-\beta_i\circ \bfF_{\delta_i}\\
	&=\alpha\circ \bfF_{\gamma}-\beta\circ \bfF_{\delta_i} +\gamma-\delta_i\\
	&=\alpha\circ \bfF_{\gamma} - \beta_i\circ \bfF_{\gamma-(c_i(g)-b)\theta_i} + \gamma-\delta_i &(\text{by Eq. \eqref{eq:diff-on-Q-C}})\\
	&=\alpha\circ \bfF_{\gamma}-\beta_i\circ \bfF_{\gamma} + \gamma-\delta_i &(\text{$\theta$ is a period function})\\
	&=(\alpha-\beta_i)\circ \bfF_{\gamma} + \gamma-\delta_i \\
	&=(c_i(h)+a)\theta_i\circ \bfF_{\gamma} +(c_i(g)+b)\theta_i\\ 
	&=(c_i(h)+a)\theta_i +  (c_i(g)+b)\theta_i & (\text{$\theta_i$ is constant on trajectories of $\bfF$}) \\
	&=(c_i(h)+c_i(g) + (a+b))\theta_i,
\end{align*}
Then, from the last formula, we obtain
\begin{align*}
[c_1(h\circ g), c_2(h\circ g), \ldots, c_n(h\circ g)] &= [c_1(h)+c_1(g), c_2(h)+c_2(g), \ldots, c_n(h) + c_n(g)],
\end{align*}
which yields $\psi(h\circ g) = \psi(h)+\psi(g)$. Thus, $\psi$ is a homomorphism.

(2) We need to show that $\ker\psi = \mathcal{S}_{\id}(f,E_f)$. Let $h$ be a diffeomorphism from $ \mathcal{G}(f, E_f)$ and $\alpha$ and $\beta_i$ be shift functions for $h$ as above, $i = 1,2\ldots, n$.

Assume that $h\in \ker\psi$, i.e., $\psi(h) = [0,0,\ldots, 0]\in \bbZ^n/\Delta$. 
Then there exists $a\in\bbZ$ such that for shift functions $\alpha-a\theta$ on $M\setminus E_f$ and $\beta_i$ on $Q_i$, Eq.  \eqref{eq:c_i(h)-per} holds, i.e., $\alpha-\beta_i = a\theta_i$, or equivalently $\alpha-a\theta = \beta_i$,  on $Q_i\setminus C_i$ for each $i = 1,2,\ldots, n$.
Since $\beta_i$ be a smooth function on $Q_i$, then a function $\alpha-a\theta$ can be extended to a smooth function $\tilde{\alpha}:M\to\bbR$ such that $h = \bfF_{\tilde{\alpha}}$. By case (4) of Proposition \ref{prop:shift-functions-exist}, a diffeomorphism  $h$ belongs to $\mathcal{S}_{\id}(f,E_f)$.

Suppose $h$ belongs to $\mathcal{S}_{\id}(f,E_f)$. Then, by case (4) of Proposition \ref{prop:shift-functions-exist}, there exist a unique function $\alpha:M\to \bbR$ such that $h = \bfF_{\alpha}$ on $M$.
 Note that $\alpha|_{M\setminus E_f} = \alpha'-a\theta$ for some smooth function $\alpha':M\setminus E_f\to \bbR$ and some $a\in \bbZ$. Then there exists a unique $c_i(h)\in\bbZ$ such that $\alpha' -a\theta-\beta_i = c_i(h)\theta_i$ on $Q_i\setminus C_i$.
Since $\alpha|_{M\setminus E_f} = \alpha'-a\theta$ and $\beta_i$ is $C^{\infty}$-functions on $Q_i$, it follows from Lemma \ref{lm:period} that there exists a sequence $\{x_{\ell}\}\subset Q_i\setminus C_i$ converging to some $ z\in C_i$ and such that $\lim_{\ell\to \infty}\theta_i(x_{\ell}) = \infty$. Therefore, $c_i(h)$ must be $0$ for all $i = 1,2,\ldots, n,$ which means that $\psi(h) = [0,0,\ldots, 0]\in \bbZ^n/\Delta$, i.e., $h\in \ker \psi.$   
\end{proof}
By Lemma \ref{lm:zeta-alpha-not-unique}, we have $\mathcal{G}(f,E_f)/\ker \psi\cong \bbZ^n/\Delta$. 
From the fact that 
$\bbZ^n/\Delta$ is isomorphic to $\bbZ^{n-1}$, together with
  $\mathcal{G}_{\id}(f, E_f) = \mathcal{S}_{\id}(f,E_f)$, we finally  obtain
$$
\pi_0\mathcal{G}(f, E_f) =  \mathcal{G}(f, E_f)/\mathcal{G}_{\id}(f, E_f)  = \mathcal{G}(f, E_f)/\mathcal{S}_{\id}(f,E_f) = \mathcal{G}(f, E_f)/\ker\psi\cong \bbZ^n/\Delta\cong \bbZ^{n-1}.
$$
This ends the proof of Proposition \ref{prop:pi0-G}.\qed

\section{Proof of Theorem \ref{thm:main-theorem}}\label{sec:proof-main}
Let $M$ be a smooth, connected, compact, and \textr{orientable} surface, and $f$ be a function from $\F$ with the set of critical points $\Crit$. Let also $E_f$ be a set of extremal circles of $f$.  
If $E_f=\varnothing$, then Theorem   \ref{thm:main-theorem} in this case 
is the statement of 
Corollary \ref{cor:hom-type-SidE}.

Assume that $|E_f| = n\geq 1$. Then by Lemma \ref{lm:weak-Sid}, $\pi_q\mathcal{S}_{\id}(f) = 0$ for $q\geq 2$ and $\pi_1\mathcal{S}_{\id}(f)$ is an abelian group for which there is a short exact sequence \eqref{eq:non-triv-subseq}. By Proposition \ref{prop:pi0-G}, $\pi_0\mathcal{G}(f,E_f)$ is a free abelian group of rank $n$, if $f$ has property (T), or $n-1$ otherwise. 
Then sequence  \eqref{eq:non-triv-subseq} always splits, and  from rank argument, we obtain
$$
\pi_1\mathcal{S}_{\id}(f)\cong \begin{cases}
	0, & \text{if }\pi_0\mathcal{G}(f, E_f)\cong \bbZ^n\\
	\bbZ, & \text{if }\pi_0\mathcal{G}(f, E_f)\cong \bbZ^{n-1}
\end{cases}
$$
Thus, $\mathcal{S}_{\id}(f)$ is {\it weakly contractible} if $f$ has property (T)  or is {\it weakly homotopy equivalent} to $S^1$ otherwise. 
By Lemma \ref{lm:CW}, $\Sid$ has the homotopy type of a CW complex. Then, by Whitehead theorem, weak homotopy equivalences obtained above are  {homotopy equivalences}.

\section{Proof of Proposition \ref{prop:main}}\label{sec:proof-main-prop}
The aim of this section is to give the proof of Proposition \ref{prop:main}.  
Let $M$ be a smooth, compact, connected, and \textr{orientable} surface, let $f$ be a function from $\F$ such that $\Sid$ is homotopy equivalent to $S^1$. By Theorem \ref{thm:main-theorem}, a function $f$ has no saddles and all isolated local extremes of $f$ are non-\textr{degenerate} (Morse). 

Throughout this section we will always assume that $\CritC = \{C_1,C_2,\ldots, C_n\}$ for some $n\geq 1$, since Proposition \ref{prop:main} is known when $|\CritC| =0$, see \cite[Remark 2.5.2.]{Maksymenko:TopAppl:2018}. Let also $E_f\subset \CritC$ be a set of extremal circles of $f$.

The following two lemmas will be needed for our proof.
\begin{lemma}\label{lm:subsuerace-N} 
	Let $N$ be a connected component of $ M\setminus \bigcup_{i = 1}^nC_i$. Then
	\begin{enumerate}
		\item $\overline{N}$ is diffeomorphic to either a cylinder or a $2$-disk. Therefore, $M$ is obtained by attaching together cylinders or/and $2$-disks along their boundary components, which are critical circles of $f$,
		\item $f|_{\overline{N}}$ has an isolated critical point {\rm(}being always a non-\textr{degenerate} local extremum and unique{\rm)} iff $\overline{N}$ is a $2$-disk. 
	\end{enumerate}
\end{lemma}
\begin{proof}
	Let $C_i\in\CritC$ be a critical circle of $f$. Denote by $Q_i$ a foliated neighborhood of $C_i$ being a cylinder. So $Q_i\setminus C_i$ has no critical points of $f$. Let $N'$ be a connected component $M\setminus \bigcup_{i = 1}^n Q_i$, which is an \textr{orientable}, compact surface with the boundary and such that $N'\subset N.$ 
	
	The restriction $g = f|_{N'}:N'\to P$ is a Morse function without saddles. 
	Hence, by Morse equalities, we have $\chi(N') = |\Sigma_{g}|\geq 0$, where $\Sigma_g$ is the set of critical points of $g$.
	Note that  $N'$ is an \textr{orientable} surface with the boundary, then $0\leq\chi(N')\leq 1$. Therefore,
	$N'$ is either a cylinder (if $\chi(N') = 0$, $f|_{N'}$ has no isolated local extrema) or a  $2$-disk $D^2$ (if $\chi(N') = 1$, $f|_{N'}$ has  a unique  isolated local extremum). 
	Since 
	$Q_i$ is a cylinder, it follows that  $\overline{N}$ is  a cylinder (or a $2$-disk) if $N'$ is a cylinder (or a $2$-disk).
\end{proof} 

The following lemma concerns functions on the torus.
{
\begin{lemma}\label{lm:torus-f}
	Let $f:T^2\to P$ be a function from the class $\mathcal{F}$ without isolated critical points and with $|E_f|\geq 2$. Let $L$ and $L'$ be two extremal circles of $f$ such that $L$ and $L'$ are connected components of the boundary of some cylinder $Q\subset T^2$ whose interior does not  contain extremal circles of $f$. Then either  $L$ is maximal and $L'$ is minimal extremal circles, or $L$ is minimal and $L'$ is maximal extremal circles.
\end{lemma} 
}
\begin{proof}
	This result follows from the properties of the gradient vector field of $f$, see \textsection\ref{ssec:class-F}. We left details to the reader.
\end{proof}

\subsection{Proof of Proposition \ref{prop:main}}
As we mentioned earlier, the result is known if $|\CritC| = \varnothing$, see \cite[Remark 2.5.2.]{Maksymenko:TopAppl:2018}.
 If $|\CritC|\geq 1$, then 
 (1) and (2) of Proposition \ref{prop:main} are corollaries of Lemma \ref{lm:subsuerace-N}.

(3)  Let us discuss some special cases. The case  $E_f = \varnothing$ is only possible if $f:T^2\to P$ is a not null-homotopic circle-valued function. It is easy to prove using properties of the gradient vector field of $f$ that the case $|E_f| = 1$ is impossible.
If $f:T^2\to P$ is null-homotopic, then from compactness of $T^2$ follows that $|E_f|\geq 2$.
It remains to show that $|E_f|= n \geq 2$ is not odd.

By Lemma \ref{lm:subsuerace-N}, $T^2$ is obtained by attaching cylinders bounded by critical circles of $f$ along their boundaries. 
Consider a courser partition of $T^2$ into cylinders $\mathcal{Q} = \{Q\subset T^2\,| \partial Q\subset E_f\}$ bounded by extremal circles $\{L\,|\, L\in E_f\}\subset \CritC$ of $f$ such in Lemma \ref{lm:torus-f}.  It is easy to see that $|\mathcal{Q}| = |E_f|$ and elements of $E_f $ and $\mathcal{Q}$ can be cyclically enumerated, i.e., $L_i = L_{i\;\mathrm{mod}\,n}$,
$Q_i = {Q}_{i\;\mathrm{mod}\,n}$, and $Q_i$ is bounded by $L_i$ and $L_{i+1},$ $i = 0,1,\ldots,n-1$.

From Lemma \ref{lm:torus-f}, in the cyclical order as above, the maximal and minimal circles in $E_f$ alternate, i.e., if $L_i$ is minimal (maximal) then $L_{i+1}$ is maximal (minimal). Assume that $|E_f| = n$ is odd and $L_0$ is maximal (minimal). Then then $L_n$ is minimal (maximal), which contradicts the fact that $L_0 = L_n$. Therefore $|E_f|$ is always even. 
\qed


\end{document}